\documentclass[numbers,webpdf]{ima-authoring-template}%

\graphicspath{{Fig/}}

\theoremstyle{plain}
\newtheorem{theorem}{Theorem}
\newtheorem{lemma}{Lemma}
\newtheorem{remark}{Remark}%

\numberwithin{equation}{section}

\begin{document}

\DOI{DOI HERE}
\copyrightyear{2021}
\vol{00}
\pubyear{2021}
\access{Advance Access Publication Date: Day Month Year}
\appnotes{Paper}
\copyrightstatement{Published by Oxford University Press on behalf of the Institute of Mathematics and its Applications. All rights reserved.}
\firstpage{1}


\title[]{Tuning Butterworth filter's parameters in SPECT reconstructions via kernel-based Bayesian optimization with a no-reference image evaluation metric}

\author{Luca Pastrello
\address{\orgdiv{Department of Mathematics\lq\lq Tullio Levi-Civita\rq\rq}, \orgname{University of Padova}, \orgaddress{\street{Via Trieste 63}, \postcode{35121}, \country{Italy}}}}
\author{Diego Cecchin\ORCID{0000-0001-7956-1924}
\address{\orgdiv{Nuclear Medicine Unit, Department of Medicine – DIMED and Padova Neuroscience Center (PNC)},\\ \orgname{University of Padova}, \orgaddress{\street{Via Giustiniani 2}, \postcode{35128}, \country{Italy}}}}
\author{Gabriele Santin\ORCID{0000-0001-6959-1070}
\address{\orgdiv{Department of Environmental Sciences, Informatics and Statistics}, \orgname{Ca’ Foscari University of Venice}, \orgaddress{\street{Via Torino 155}, \postcode{30172}, \country{Italy}}}}
\author{Francesco Marchetti*\ORCID{0000-0003-1087-7589}
\address{\orgdiv{Department of Mathematics \lq\lq Tullio Levi-Civita\rq\rq and Padova Neuroscience Center (PNC)},\\ \orgname{University of Padova}, \orgaddress{\street{Via Trieste 63}, \postcode{35121}, \country{Italy}}}}

\authormark{Pastrello et al.}

\corresp[*]{Corresponding author: \href{email:francesco.marchetti@unipd.it}{francesco.marchetti@unipd.it}}

\received{Date}{0}{Year}
\revised{Date}{0}{Year}
\accepted{Date}{0}{Year}


\abstract{In Single Photon Emission Computed Tomography (SPECT), the image reconstruction process involves many tunable parameters that have a significant impact on the quality of the resulting clinical images. Traditional image quality evaluation often relies on expert judgment and full-reference metrics such as MSE and SSIM. However, these approaches are limited by their subjectivity or the need for a ground-truth image. In this paper, we investigate the usage of a no-reference image quality assessment method tailored for SPECT imaging, employing the Perception-based Image QUality Evaluator (PIQUE) score. Precisely, we propose a novel application of PIQUE in evaluating SPECT images reconstructed via filtered backprojection using a parameter-dependent Butterworth filter. For the optimization of filter's parameters, we adopt a kernel-based Bayesian optimization framework grounded in reproducing kernel Hilbert space theory, highlighting the connections to recent greedy approximation techniques. Experimental results in a concrete clinical setting for SPECT imaging show the potential of this optimization approach for an objective and quantitative assessment of image quality, without requiring a reference image.}


\maketitle

\section{Introduction}

In medical imaging, clinical images that are evaluated by the expert physician are the result of a reconstruction process that naturally encompasses the tuning and selection of many (hyper)parameters. This is also true in  Single Photon Emission Computed Tomography (SPECT), which relies on the detection of emissions from radiopharmaceuticals introduced into the body, and provides functional and metabolic information. Together with Positron Emission Tomography (PET), they represent two fundamental imaging techniques in nuclear medicine, sharing similar image reconstruction algorithms to a certain extent \cite{Bailey05,Cherry12}.

Differently with respect to PET, where the annihilation of positron-emitting radionuclides produces pairs of gamma-photons that are detected in coincidence, SPECT imaging uses single-photon-emitting radionuclides like technetium-99m, and has generally lower spatial resolution than PET. Nevertheless, its cost-effectiveness and availability make SPECT a valuable alternative to PET imaging in many medical areas such as oncology and neurology, where it plays a crucial role in disease diagnostic and monitoring \cite{Rahmim08,Zanzonico04}.  Its clinical utility is further enhance by hybrid imaging approaches (e.g. SPECT/CT), where functional and anatomical information are combined \cite{Hutton19}.

Very often, the evaluation of the goodness of a reconstructed clinical image is carried out by the physicians, who consequently set up the parameters' 
framework in light of their expertise and of some perceived image quality, and by taking into account possible guidelines offered by relevant healthcare 
companies. However, this approach is very likely to be strongly user-dependent, and therefore a \textit{quantitative} approach in image quality assessment 
should be preferable in order to favor some shared \textit{consensus}.

In this direction, full-reference metrics such as, e.g., the Mean Squared Error (MSE) and the Structural Similarity Index (SSIM) are commonly used in medical image evaluation. They are capable of capturing many types of degradations and distortions that may take place along the processing step (e.g., during acquisition, compression, etc), but their major drawback consists in their need of a reference image, a \textit{ground truth}, to be used in their computation. Unfortunately, this severely limits their concrete applicability in the clinical practice, as it is usually very hard to get a target image, which may require common standardization, preprocessing and registration strategies, as it could be also produced by using different medical devices. 

Therefore, no-reference metrics may play a crucial role, as they are meant to not require the presence of any target image, making them particularly versatile for many applications. In \cite{mittal2012}, Mean Subtracted Contrast Normalized (MSCN) coefficients were introduced as a tool to reach a satisfactory no-reference image quality assessment. Based on this, the MSCN coefficients were later utilized in \cite{venkatanath2015} to define a comprehensive no-reference metric known as Perception-based Image QUality Evaluator (PIQUE). This metric uses perception-based features to assess image quality without requiring a reference image, and its usage in the PET imaging framework was introduced and discussed in \cite{shigeaki2024}, where its effectiveness in assessing the quality of reconstructed medical images was highlighted. 

From these considerations, the contribution of this paper is the following. First, we consider for the first time a no-reference image quality assessment approach in the context of SPECT imaging. Specifically, we use the PIQUE score to assess the quality of filtered back-projected images, where the considered Butterworth filter depends on some parameters that need to be tuned. For this purpose, here we consider a kernel-based Bayesian Optimization (BO) approach, which is a well established gradient-free optimization technique that is grounded in the theory of reproducing kernel Hilbert spaces, as we describe in Section \ref{sec:kernel} providing an original link to novel greedy approaches existing in the framework of kernel-based approximation. After a detailed introduction to the considered PIQUE score in Section \ref{sec:pique}, we present the results obtained in a concrete experimental SPECT imaging setting in Section \ref{sec:results}.

\section{Kernel-based Bayesian optimization}\label{sec:kernel}

\subsection{Approximation by positive definite kernels}

For the following overview on classical topics of kernel-based approximation theory, we refer e.g. to the well-known books \cite{Fasshauer07,Fasshauer15,Wendland05}.

In the classical RKHS approximation theory framework, the aim is to approximate an unknown function $f: \Omega \longrightarrow \mathbb{R}$, 
$\Omega\subset\mathbb{R}^d$, by means of a function $\hat{f}^{(n)}\in\mathcal{K}_n=\mathrm{span}\{\kappa(\cdot,\mathbf{x}_i),\;\mathbf{x}_i\in X_n\}$, where 
$\kappa:\Omega\times\Omega\to\mathbb{R}$ is a (strictly) positive definite kernel and $X_n=\{\mathbf{x}_1,\dots,\mathbf{x}_n\}\subset\Omega$ is a set of data 
sites, where knowledge on the unknown function is assumed. Usually, the kernel $\kappa$ also depends on a shape parameter $\varepsilon>0$ (e.g., the well-known 
Gaussian kernel $\kappa_G(\mathbf{x},\mathbf{y})=e^{-(\varepsilon\lVert \mathbf{x}-\mathbf{y}\lVert)^2}$), which plays the role of a \textit{hyperparameter} of 
the approximation process. Generalization to multiple hyperparameters, and resulting anisotropic kernels, have been discussed in literature too \cite{Wenzel24}. 
By denoting $Z_n=(X_n,f(X_n))$, the approximation problem can be translated into the optimization task \vspace{-0.1 cm}
\begin{equation}\label{eq:mini}
	\hat{f}^{(n)}=\hat{f}_R^{(n)}=\textrm{argmin}_{f'\in\mathcal{H}_\kappa} R_{Z_n}(f'), \vspace{-0.1 cm}
\end{equation}
where the completion $\mathcal{H}_\kappa=\overline{\textrm{span}\{ \kappa(\cdot,\mathbf{x}), ~ \mathbf{x}\in\Omega\}}$ is the RKHS of $\kappa$ on $\Omega$, and $R_{Z_n}$ is an appropriate loss function to be minimized. 
Particularly, we consider the square loss $$R_{Z_n}(f')=\sum_{i=1}^n (f(\mathbf{x}_i)-f'(\mathbf{x}_i))^2.$$
In this case, if $\kappa$ is strictly positive definite, we observe that $\hat{f}^{(n)}$ in \eqref{eq:mini} is the kernel-based interpolant of $f$ at $X_n$, i.e., $\hat{f}^{(n)}(\mathbf{x})=\sum_{i=1}^n c_i \kappa(\mathbf{x},\mathbf{x}_i)$ where the coefficients $c_1,\dots,c_n$ are uniquely determined by the interpolation conditions. 
Such interpolatory framework has been extensively studied in the literature, and various error bounds for this reconstruction process have been designed. 
For example, letting $\kappa_\mathbf{x}(\cdot)=\kappa(\cdot,\mathbf{x})$ and 
$\kappa_\mathbf{x}({X_n})=(\kappa_\mathbf{x}(\mathbf{x}_1),\dots,\kappa_\mathbf{x}(\mathbf{x}_n))$, a well-known pointwise error bound involves the so-called 
\textit{power function} \vspace{-0.1 cm} $$
P_{\kappa,X_n}(x)=\left(\kappa_\mathbf{x}(x) - \kappa_\mathbf{x}({{X}_n})\cdot(\mathsf{K}^{-1} \kappa_\mathbf{x}({{X}_n}))\right)^{1/2},
$$ 
where 
$\mathsf{K}_{ij}=\kappa(\mathbf{x}_i,\mathbf{x}_j)$ is the kernel matrix, and reads as follows \cite[Theorem 14.2, p. 117]{Fasshauer07} for $f\in\mathcal 
H_\kappa$: \vspace{-0.1 cm}
\begin{equation}\label{eq:the_error}
	|f(\mathbf{x})-\hat{f}^{(n)}(\mathbf{x})|\le P_{\kappa,X_n}(\mathbf{x})\lVert f\lVert_{\mathcal{H}_\kappa}, \quad \mathbf{x}\in\Omega. \vspace{-0.1 cm}
\end{equation}
Moreover, by defining the \textit{fill distance} $h_{{X}_n,\Omega} =  \sup_{ \mathbf{x} \in \Omega} {\min_{ \mathbf{x}_i  \in {{X}_n}} \lVert \mathbf{x} - 
\mathbf{x}_i \lVert_2 }$, assuming $\Omega$ to be bounded and satisfying an interior cone condition, and $\kappa \in C^{2s}\left(\Omega \times \Omega\right)$, 
we then obtain for \textit{small enough} $h_{{X}_n,\Omega}<h_0$, $h_0>0$, the further bound \vspace{-0.1 cm}
\begin{equation}\label{eq:conv_rate}
	|f\left({\mathbf{x}}\right)-\hat{f}^{(n)}(\mathbf{x})| \leq C_\kappa(\mathbf{x}) h^s_{{X}_n,\Omega}  \lVert f \lVert_{\mathcal{H}_\kappa},\vspace{-0.1 cm}
\end{equation}
where the factor $C_\kappa(\mathbf{x})$ depends on the maximum of kernel derivatives of degree $2s$ in a neighborhood of $\mathbf{x}\in\Omega$ \cite[Section 14.5]{Fasshauer07}. 
On the one hand, the theoretically achievable convergence rate is therefore influenced by the fill distance and the smoothness of the kernel. On the other hand, two terms play an important role in affecting the conditioning of the interpolation process. 

First, the \textit{separation distance} $q_{X} = \frac{1}{2} \min_{  i \neq j} \lVert {\mathbf{x}}_i - {\mathbf{x}}_j \lVert_2$. The interpolation process gets more ill-conditioned as the separation distance becomes smaller, which is what usually happens in practice when increasing the number of data sites in order to reduce the value of the fill distance. Moreover, choosing a regular kernel theoretically leads to a faster convergence rate, but the more regular the kernel, the larger the ill-conditioning in the interpolation process. These are often referred to as \textit{trade-off} principles in the relevant literature (see e.g. \cite[Section 16]{Fasshauer07} for an overview). 

Second, the value of the shape parameter $\varepsilon$ of the kernel $\kappa$. A large value produces very localized basis functions that lead to a well-conditioned but likely inaccurate approximation scheme, while by lowering such value we may obtain a more accurate reconstruction at the price of an ill-conditioned setting. Consequently, accurately tuning kernel's hyperparameters is an important issue analyzed in many less and more recent works, e.g. \cite{Cavoretto18,Chen22,Fornberg07,Marchetti24,Scheuerer11}. Multigrid approaches with different lengthscale parameters in place of a unique shape parameter have been designed too \cite{Wendland10}.

In case the approximant $\hat{f}^{(n)}$ can not be found in a stable way from a numerical point of view, and/or the values $f(X_n)$ are affected by noise, a 
possible approach consists of regularizing the approximation problem \eqref{eq:mini} via a Tikhonov-type formulation. In this framework, to obtain error 
bounds which share the same spirit of \eqref{eq:conv_rate}, \textit{sampling inequalities} represent a fundamental tool in the literature 
\cite{Narcowich05,Rieger09,Wendland05a}. Their idea consists of loosing the interpolatory constraints and to admit some reasonable \textit{small} error at the 
data sites $X$. 

Another strategy that has been considered to avoid instability issues in kernel methods, and possibly increase their accuracy, consists in constructing the 
model at selected nodes $X_m$ extracted from the original set of data sites $X_n$, with $m\ll n$. In this manner, it is possible to obtain a \textit{sparse} 
approximation of the RKHS approximant, leading to an evident advantage in terms of model complexity and evaluation speed. Greedy kernel methods represent a well 
established approach for selecting a meaningful subset $X_m \subset X_n$, by following an iterative procedure that starts from an empty set $X_0 = \emptyset$ 
\cite{Wenzel23,Wirtz15}. 
At the $m$-th step, the set $X_m$ is defined as $X_{m} = X_{m-1} \cup \{{x}_m\}$, where ${x}_m$ is selected so that 
\begin{equation}\label{eq:optimus}
	{\mathbf{x}}_{m} = \textrm{argmax}_{{\mathbf{x}} \in X_n \setminus X_{m-1}} \eta^{(m)}({\mathbf{x}}),
\end{equation}
where $\eta^{(m)}: \Omega \longrightarrow \mathbb{R}$ is an error indicator. In literature, the most frequent error criteria are based on either the power 
function (see \eqref{eq:the_error}) or the residual $r_{X_m} \equiv f - \hat{f}^{(m)}$ :
\begin{enumerate}
	\item $P$-greedy: $\eta_P^{(m)}(\boldsymbol{\mathbf{x}}) = P_{\kappa,X_m}(\boldsymbol{\mathbf{x}})$, promotes the selection of quasi-uniform data sites \cite{DeMarchi05}. A thorough analysis of the convergence rates can be found in \cite{Santin16, Wenzel21}.
	\item $f$-greedy: $\eta_f^{(m)}(\mathbf{x}) = |r_{X_m}(\mathbf{x})|$, enforces the choice of data sites according to the structure of the underlying function \cite{Schaback00}.
\end{enumerate}
Very recently, $P$-greedy, $f$-greedy and combinations of the two criteria were unified within the framework of $\beta$-greedy algorithms \cite{Wenzel23}. These 
adaptive approaches for function approximation have been gaining a consistent interest in the relevant literature during the last years, due to their 
flexibility compared to non-adaptive interpolation and theoretically certified accuracy \cite{Santin24}. Also, their usage is facilitate by a 
well established numerical software \cite{Wirtz13}.

\subsection{Bayesian optimization approaches}

In addition to function approximation, greedy strategies built upon kernel models are considered for addressing optimization tasks, where such approaches fall within the framework of kernel-based Bayesian Optimization (BO), as it is referred to as in the relevant literature. In this case, the aim consists in estimating the maximum point ${\mathbf{x}}^{*}$ of the unknown function $f:\Omega\longrightarrow\mathbb{R}$, i.e., 
\begin{align*}
	{\mathbf{x}}^{*}= \textrm{argmax}_{\boldsymbol{\mathbf{x}} \in \Omega}\:f(\boldsymbol{\mathbf{x}}).
\end{align*}
The approach carried out by BO is model-based, meaning that the searching procedure is guided with the help of a surrogate model, and gradient-free, as no knowledge on derivatives of $f$, nor approximations, are considered. We remark that the function $f$ might be non-concave, thus it can admit multiple local maxima. Moreover, evaluating $f$ during the searching procedure may be expensive in many applications, and therefore BO needs to consider the issue of finding a satisfactory balance between \textit{exploration} (i.e., gathering information on $f$ in $\Omega$) and \textit{exploitation} (i.e., using the constructed surrogate model to estimate ${\mathbf{x}}^{*}$). 

We briefly present BO while staying on track with the discussion about the greedy approach. Given a \textit{large} set $X_n\subset\Omega$ of data sites, starting from a \textit{small} subset $X_0\subset X_n$ (possibly the empty set), at step $m$ the set $X_m$ is again defined as $X_{m} = X_{m-1} \cup \{{x}_m\}$ where ${\mathbf{x}}_m$ is selected by taking (cf. \eqref{eq:optimus})
\begin{equation}\label{eq:rule_bo}
	\eta^{(m)}({\mathbf{x}})=\hat{f}^{(m)}(\mathbf{x})+\beta_m P^2_{\kappa,X_m}(\mathbf{x}),
\end{equation}
where $\hat{f}^{(m)}(\mathbf{x})$ is the kernel approximant constructed at $X_m$ and $\beta_m>0$ is a parameter that rules the trade-off between exploration and 
exploitation. We point out that interpreting the power function in the shape of a confidence interval, as done in BO, is linked to the statistical perspective 
of data interpolation \cite{kanagawa2018} also provided by the \textit{kriging} method \cite{Scheuerer13}. Moreover, approaches where not only a single node, 
but a whole batch of data sites is added at each step have been considered in literature, including the noisy setting.

Instead of looking at the approximation error, in this context the effectiveness of the selection strategy is usually measured by considering the cumulative regret $\textrm{c-reg}_{m}$ and the simple regret $\textrm{s-reg}_{m}$, which are defined as
\begin{equation*}
	\textrm{c-reg}_{m}=\sum \limits_{t=1}^{m} (f(\boldsymbol{\mathbf{x}}^{*})-f(\boldsymbol{\mathbf{x}}_{t})),\quad
	\textrm{s-reg}_{m}=f(\boldsymbol{\mathbf{x}}^{*})-\max_{\mathbf{x}\in X_m}f(\mathbf{x}).
\end{equation*}
While $\textrm{c-reg}_{m}$ quantifies an overall cost that depends on the computed queries, the actual convergence of the iterative scheme is expressed in terms of $\textrm{s-reg}_{m}$, and it is influenced by the choice of $\beta_m$. Note that the trade-off between exploitation and exploration realized by $\beta_m$ needs to take into account the setting of the chosen application, and the corresponding cost of exploratory queries. A careful design of this trade-off strategy is a crucial issue for this optimization approach, which has been observed and discussed in many works \cite{Bull11,Lyu19,Srinivas10}. In fact, although various theoretical bounds were obtained in the relevant literature, we point out that precise convergence estimates for $\textrm{c-reg}_{m}$ and $\textrm{s-reg}_{m}$ that fully exploit the choice of weights and the regularity of the kernel $\kappa$ still represent an open problem in the relevant literature \cite{Vakili22}, and an under-developed research line if compared to the greedy approximation framework.

\begin{remark}\label{rem:rob}
	The initialization method of the greedy optimization, that is, the definition of $X_0$, has been discussed e.g. in \cite{Lyu19}. It turns out that a robust initialization shares a similar spirit with known results from approximation theory such as \eqref{eq:the_error}, meaning that a satisfactory covering of the domain (i.e., small fill distance) is desirable.
\end{remark}

\section{No-reference perception-based image quality evaluator}\label{sec:pique}

For the sake of completeness, in this section we outline the main steps required to define the Perception-based Image QUality Evaluator (PIQUE) no-reference metric, which we will consider in our image reconstruction tasks. We refer to the seminal paper \cite{venkatanath2015} for more details.

First, some preprocessing steps are carried out before actually computing the metric on an image, starting from local mean removal and divisive normalization. 
This step is used to distinguish the blocks in the image into uniform (U) and spatially active (SA). This distinction is crucial because only the blocks are analyzed further, given that humans predominantly focus on these regions, according to the HSV model.

In practice, let $I(i, j)$ be the image luminance, where $i = 1,\dots,M$ and $j=1,\dots, N$ are spatial indices, and $M,\;N$ are the height and width of the image, respectively. We then define the Mean Subtracted Contrast Normalized (MSCN) coefficients
\begin{align}\label{MSCN}
	\hat{I}(i, j)= \frac{I(i, j) -\mu(i, j)}{\sigma(i, j)+C},
\end{align}
where $C$ is a constant useful in preventing instability issues, and
\begin{align*}
	\mu(i, j)= \sum_{k=-K}^{K}\sum_{l=-L}^{L}w_{k, l}I_{k, l}(i, j),\quad \sigma(i, j)= \sqrt{ \sum_{k=-K}^{K}\sum_{l=-L}^{L}w_{k, l}(I_{k, l}(i, j)- \mu(i, j))^{2}},
\end{align*}
where $w=\{ w_{k,l}\:|\: k=-K, \dots, K,\;l=-L,\dots, L \}$ is a 2D circularly symmetric Gaussian weighting function sampled out to three standard deviations $(K=L=3)$ and rescaled to unit volume. The two scalar fields $\mu$ and $\sigma $ are called the local mean field and local deviation field, respectively.

After computing the MSCN coefficients, the image is segmented into non-overlapping blocks $B_{k}$ of size $n \times n$, $k=1,\dots,N_B$, leaving out those at the image boundaries. The  MSCN coefficients are utilized to label each block as a uniform (U) or as a spatially active (SA) block. The labeling criterion is
\begin{align}
	B_{k}=
	\begin{cases}
		\textrm{U},   &  \nu_{k} <T_{U},\\
		\textrm{SA},   &  \nu_{k}\geq T_{U},
	\end{cases}
\end{align}
where $\nu_{k}$ is the variance of the MSCN coefficients $\hat{I}(i, j)$ of a given block $B_{k}$, that is
\begin{align}
	\nu_{k}= \frac{1}{n^2}\sum_{i=1}^{n}\sum_{j=1}^{n} \big( \hat{I}(i, j)-\frac{1}{n^2}\sum_{i=1}^{n}\sum_{j=1}^{n}\hat{I}(i, j) \big)^{2}, \quad k \in 1, 2, \dots, N_{B},
\end{align}
and $T_{U}$ is a threshold value, which can be expressed as a percentage as we are dealing with normalized coefficients.

The second preprocessing step consists of analyzing each SA block and assign a score with respect to two types of distortion criteria, namely, \textit{noticeable distortion} criterion and \textit{additive white noise} criterion.

A block-level distortion is noticeable if at least one of its edge segments (defined as a collection of $m$ contiguous pixels in a block edge) exhibits low 
spatial activity. Precisely, each edge $L_{p}$, $p\in \{1, 2, 3, 4\}$,  is divided into segments as $a_{pq}=\{ L_{p}(x)\::\: x=q, \dots, q+(m-1)\},$ where 
$a_{qp}$ is the structuring element, and $ q \in \{1, 2, \dots, n-m\}$, denotes the segment index. We say that a segment exhibits low spatial activity if the 
standard deviation of the segment $a_{qp} $ is less than a fixed threshold $T_{STD}$, i.e.,
\begin{align}\label{condition for a distorted segment}
	\sigma_{pq}<T_{STD}.
\end{align}
Then, a block is considered distorted if at least one of its segments is linked to low spatial activity.

As far as the white noise criterion is concerned, this is inspired by the Human Visual System's (HVS) sensitivity to center-surround variations. Each block is divided into two parts, i.e., a central part $S_{cen}$ that contains the two central columns and the set of remaining columns $S_{sur}$. By defining
\begin{align*}
	\beta= \frac{|(\frac{\sigma_{cen}}{\sigma_{sur}})-\sigma_{k}|}{\max(\frac{\sigma_{cen}}{\sigma_{sur}},\sigma_{k})},
\end{align*}
where $\sigma_{cen}$, $\sigma_{sur}$ and $\sigma_{k}$ are the standard deviations of $S_{cen}$, $S_{sur}$ and of the block $B_{k}$, respectively, it was 
empirically observed that a given block can be categorized as affected by noise if it satisfies the condition
\begin{align}\label{block affected with noise condition}
	\sigma_{k}>2 \beta.
\end{align}
Now, the key insight is that the variance $\nu_{k}$ of the MSCN coefficients of a given block is a significant indicator of the amount of distortion present in that block. This can be observed e.g. in JPEG compression and blur distortion, where the variance is very low for highly distorted images when compared to images of better quality, indicating that variance is inversely proportional to the amount of distortion. On the other hand, as far as additive noise is concerned, the variance is expected to increase with increasing noise, being directly proportional to the amount of distortion. Therefore, we use $\nu_{k}$ to assign a score $D_{k}$ for each distorted block $B_{k}$ as follows
\begin{align*}
	D_{k}=
	\begin{cases}   
		1,  &\text{if} \quad \eqref{condition for a distorted segment} \quad \text{and} \quad \eqref{block affected with noise condition},\\
		\nu_{k}, & \text{if} \quad \eqref{block affected with noise condition},\\
		(1-\nu_{k}), &\text{if} \quad \eqref{condition for a distorted segment}.
	\end{cases}
\end{align*}
The PIQUE is then given by
\begin{align*}
	\textrm{PIQUE}=  \frac{(\sum_{k=1}^{N_{SA}}D_{k})+C}{N_{SA}+C},
\end{align*}
where $N_{SA}$ is the number of SA blocks in a given image, and where $C$ is a positive constant that is included to prevent numerical instability. The PIQUE 
score lies in the range from 0 to 1: the smaller the value of the index, the higher the quality of the image.

\section{SPECT image reconstruction setting}\label{sec:results}

\subsection{The StarGuide camera}

\begin{figure}[h]
	\centering
	\includegraphics[width=0.7\textwidth]{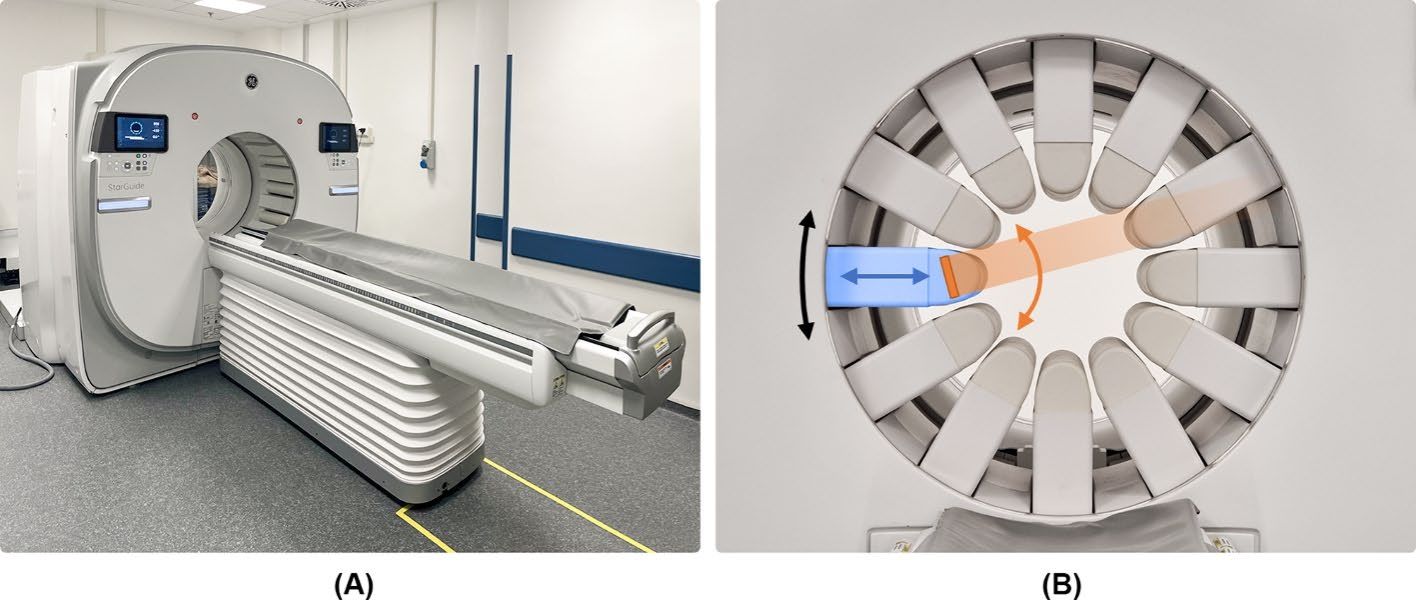}
	\caption{The StarGuide CZT SPECT/CT system. Left: an overview of the scanner. Right: the multi-detector system (from \cite{Zorz24}).}
	\label{fig:starguide}
\end{figure}

In our experiments, we make use of \textit{raw data} produced by the StarGuide SPECT/CT camera that is available at the advanced clinical unit of Nuclear Medicine, University/Hospital of Padova (see Figure \ref{fig:starguide}, left). Produced by GE HealthCare, the StarGuide is one of the only two 3D-ring general-purpose  Cadmium Zinc Telluride (CZT) SPECT/CT currently commercially available \cite{LeRouzic21}. Its multi-detector system consist of 12 columns arranged in a 3D-ring configuration over an 80 cm diameter bore. Each column consists of seven modules of 16 × 16 CZT pixelated crystals with a dual-pitched integrated parallel-hole tungsten collimator. The collimator septa are aligned with each detector pixel. The detector thickness is 7.25 mm, while the size of each of the 16 × 16 pixelated detector elements is 2.46 × 2.46 ­mm$^2$. It is important to remark that each detector column moves independently from the others, and has different degrees of movements: an automated radial
motion (in and out), a rotational motion along the gantry (from 2 to 6 steps for each bed position with a rotation range up to 25°) and a sweep motion (with a sweep range up to ± 15° in a step and shoot or continuous sweep mode) (see Figure \ref{fig:starguide}, right) \cite{Zorz24}.

\subsection{Filtered back projection}

Various image reconstruction methods have been considered in SPECT imaging. Iterative algorithms such as Maximum Likelihood Expectation Maximization (ML-EM) and derived Ordered Subset EM (OSEM) have been the standard strategies in SPECT imaging for many years because of their robustness and effectiveness \cite{Hutton11,Reader07,Zeng10}. Nevertheless, more traditional analytical approaches still represent a meaningful research line for both theoretical and applied studies, as they are involved in the initialization of iterative schemes, and explicitly deal with the physical nature of the projection data \cite{Kak01,Natterer86}.

\begin{figure}[h]
	\centering
	\begin{minipage}{0.45\textwidth}
		\centering
		\includegraphics[height=7cm]{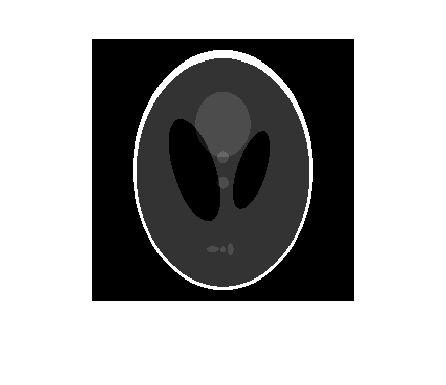} 
	\end{minipage}
	\begin{minipage}{0.45\textwidth}
		\centering
		\includegraphics[height=7cm]{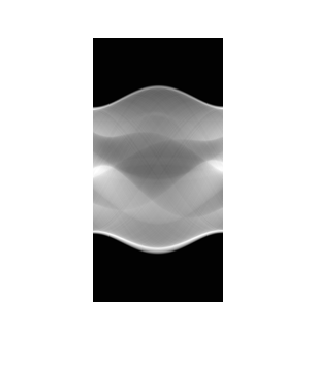} 
	\end{minipage}
	\caption{The Shepp-Logan phantom (left) and its corresponding sinogram (right) obtained through the Radon transform.}
	\label{fig:logan}
\end{figure}

The most well-known analytical approach is the so called Filtered Back Projection (FBP), where the true radiotracer distribution is estimated from measured 
projection data by using the Radon transform. We briefly recall that the projection data \( \mathcal{P}[H](\theta,s) \) at angle \( \theta \) and distance $s$ 
along the projection is given by the Radon transform of the object $H$, that is
\begin{equation*}
	\mathcal{P}[H](\theta,s) = \int_{-\infty}^{\infty} H(s \cos\theta - t \sin\theta, s \sin\theta + t \cos\theta) \, \mathrm{d}t.
\end{equation*}
The outcome of the Radon transform is usually called the \textit{sinogram} (see Figure \ref{fig:logan}).
Back-projecting the projections without filtering would result in
\begin{equation*}
	\mathcal{B}[\mathcal{P}[H]](x, y) = \int_{0}^{\pi} \mathcal{P}[H](\theta,x \cos\theta + y \sin\theta) \, \mathrm{d}\theta.
\end{equation*}

However, in this manner low-frequency components are overly emphasized, and blurring artifacts are a consequence. To avoid this issue, in FBP a combination of filters is usually applied in the frequency domain. Here, we consider the product of the high-pass \textit{ramp filter} with a low-pass  Butterworth filter, that is
\begin{equation}\label{eq:filter}
	\tau_{\rho,\omega_0}(\omega)=|\omega|\sqrt{\frac{1}{1+(\frac{\omega}{\omega_0})^{2\rho}}},
\end{equation}
where $\omega_0$ is the \textit{critical frequency} and $\rho$ is the \textit{order}. After applying $\tau_{\rho,\omega_0}$ to $\mathcal{P}[H]$ in the frequency domain, thus obtaining the filtered $\overline{\mathcal{P}[H]}$, we can thus define the filtered back-projection $\mathcal{B}_{\rho,\omega_0}[\overline{\mathcal{P}[H]}]$.

\section{Results}
\subsection{Parameters' tuning}

In the following, we accurately describe our test setting. We are considering a greedy scheme as in \eqref{eq:optimus}, where our target function $f$ consists of the opposite of the PIQUE value calculated on a filtered back-projected 3D volume $\mathcal{B}_{\rho,\omega_0}[\overline{\mathcal{P}[H]}](x,y,z)$. Precisely, we consider the mean of the PIQUE values computed for each $z$-axis slide, i.e.,
\begin{equation*}
	f(\rho,\omega_0) = -\frac{1}{Z}\sum_{z}\textrm{PIQUE}(\mathcal{B}_{\rho,\omega_0}[\overline{\mathcal{P}[H]}](x,y,z)),
\end{equation*}
where $Z$ is the number of $z$ slices, and the resulting maximization problem
\begin{equation}\label{eq:max}
	(\rho^*,\omega_0^*)= \textrm{argmax}_{\rho,\omega_0}\:f(\rho,\omega_0).
\end{equation}
To maximize $f$, we rely on the greedy rule \eqref{eq:rule_bo} and we compare three different trade-off criteria:
\begin{itemize}
	\item 
	A constant weight $\beta_m\equiv \beta$, with $\beta=\lVert f \lVert_{\mathcal{H}_\kappa}$ was analyzed in \cite{Lyu19}, and related theoretical bounds were obtained.
	\item 
	A varying increasing weight $\beta_m=\mathcal{O}(\sqrt{\log m})$ was investigated in \cite{Srinivas10}. According to the discussion carried out in the reference paper, here we choose $\beta_m=\sqrt{\log( \frac{10}{3}m^2\pi^2)}$.
	\item 
	A decreasing weight $\beta_m=\lambda^{m-1}\lVert f \lVert_{\mathcal{H}_\kappa}$, $0<\lambda<1$. We choose $\lambda=0.9$.
\end{itemize}
We highlight that the three considered criteria represent three different approaches concerning the trade-off exploration vs exploitation.

As far as the kernel is concerned, we use the well-known Matérn kernel $\kappa_M(\mathbf{x},\mathbf{y})=e^{-\varepsilon\lVert \mathbf{x}-\mathbf{y}\lVert}$, 
and we set $\varepsilon=0.1$. 

We consider three different strategies for the initialization set $X_0$ (cf. Remark \ref{rem:rob}):
\begin{itemize}
    \item 
    $X_0=\{(5,0.5)\}$, so $|X_0|=1$.
    \item 
    $X_0=\{4,8\}\times\{0.4,0.8\}$, so $|X_0|=4$.
    \item 
    $X_0=\{3,6,9\}\times\{0.3,0.6,0.9\}$, so $|X_0|=9$.
\end{itemize}
Our purpose is to avoid the calculation of too many FBP reconstructions for tuning the filter's parameters. Therefore, we limit the maximum number of observations to $20$. Consequently, the greedy algorithm will perform a total number of $20-|X_0|$ steps: depending on the initialization, we will assume more knowledge on the underlying function and reserve a restricted number of algorithm iterations, or alternatively, we will start with very limited information on the function and give more space for algorithm's exploration.

 The optimization \eqref{eq:max} is performed by searching the parameters on a fine equispaced grid
 \begin{equation*}
     \Xi=\{\rho^1,\dots,\rho^M\:|\: \rho^1=1,\;\rho^M=10\}\times \{\omega_0^1,\dots,\omega_0^M\:|\: \omega_0^1=0.1,\;\omega_0^M=1\}
 \end{equation*}
with $M=1000$. Moreover, the RKHS norm of the underlying function is estimated by computing the RKHS norm of the kernel interpolant.

The experiments were carried out in Python and PIQUE was calculated using the code available at \texttt{ https://pypi.org/project/pypiqe/} \cite{venkatanath2015}.

\subsection{Experimenting with a phantom}

In the following, we apply our framework for the reconstruction of the well-known Jaszczak SPECT/PET Phantom \cite{Jaszczak80}, with both hot and cold spheres in a hot-background (Tc-99m radiotracer ratio at least 5:1).

\begin{figure}
	\centering
		\includegraphics[width=0.99\textwidth]{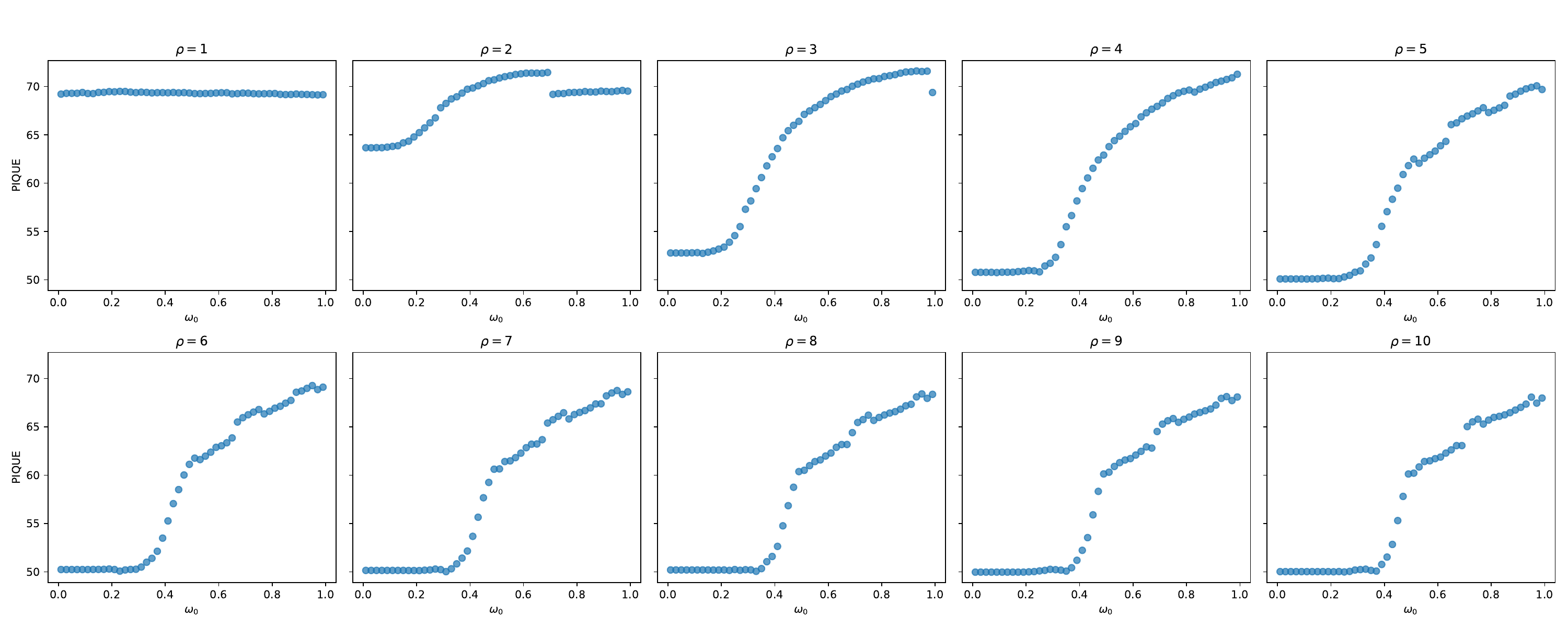}
	\caption{Values varying $\rho$ and $\omega_0$. The obtained minimum value of PIQUE is $49.966$ at $\rho=9$ and $\omega_0=0.17$.}
	\label{fig:sfere_chart}
\end{figure}

\begin{figure}
	\centering
		\includegraphics[width=0.22\textwidth]{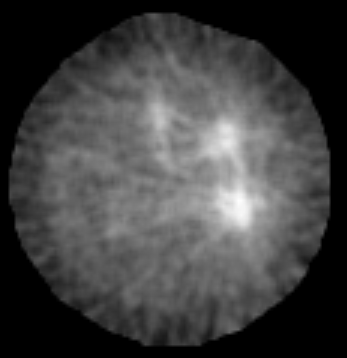}
	\includegraphics[width=0.22\textwidth]{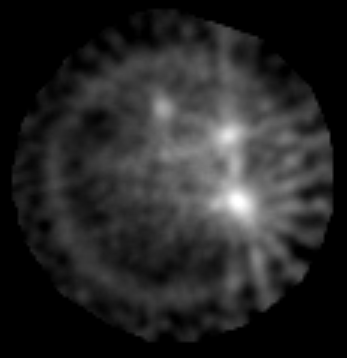}
    \caption{FBP reconstructions for a fixed z-slice of the phantom. Left: $\rho=6$, $\omega_0=0.9$. Right: $\rho=9$ and $\omega_0=0.17$ (the optimal one in Figure \ref{fig:sfere_chart}).}
	\label{fig:fbp_compare}
\end{figure}

\begin{figure}
    \centering
    \begin{minipage}[t]{0.28\textwidth}
        \centering
        \includegraphics[width=\textwidth]{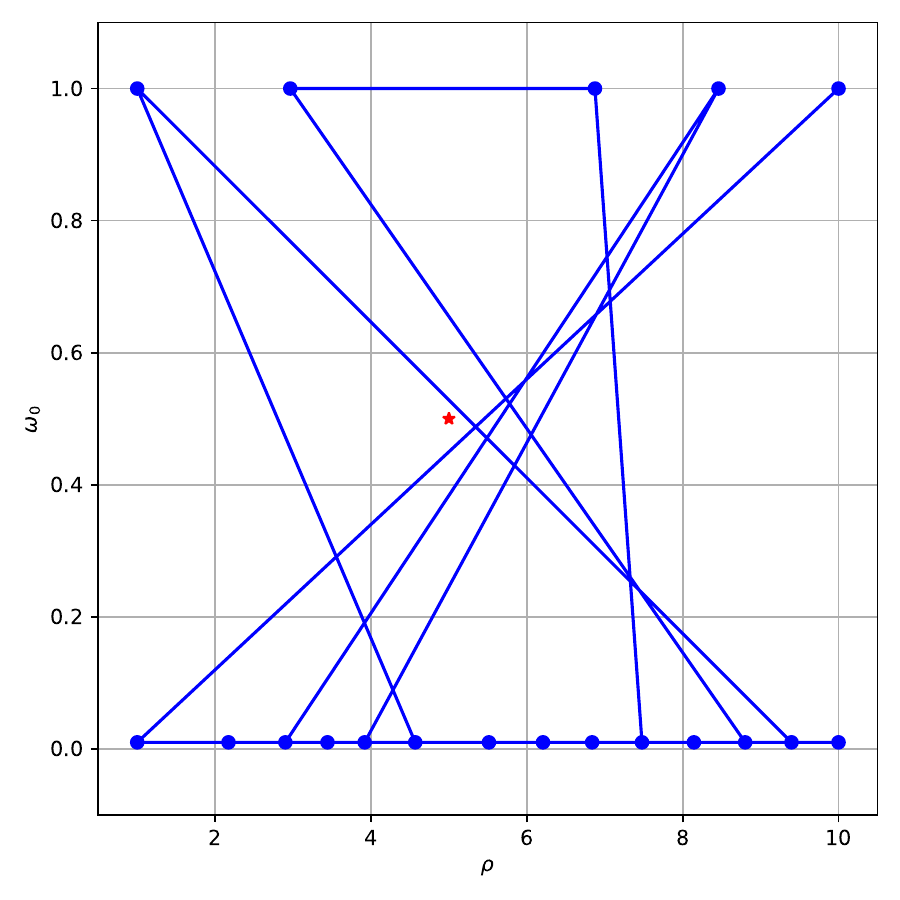}
        {\small (a) Minimum PIQUE: $49.979$.}
    \end{minipage}
    \hfill
    \begin{minipage}[t]{0.28\textwidth}
        \centering
        \includegraphics[width=\textwidth]{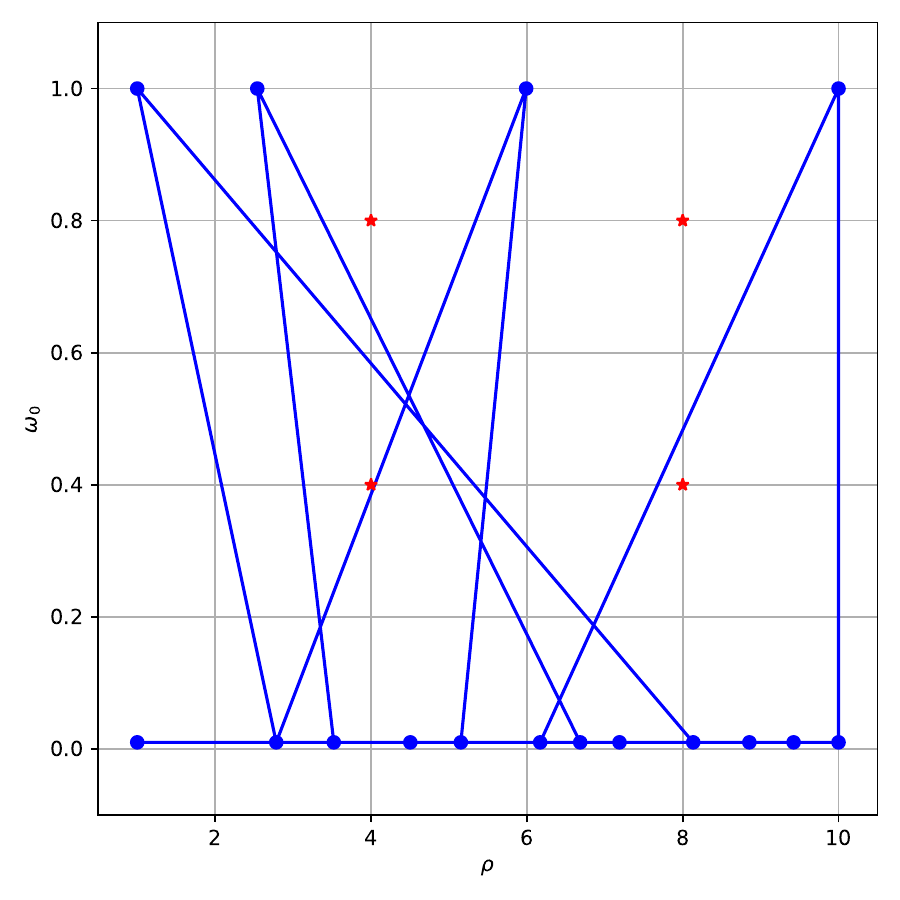}
        {\small (b) Minimum PIQUE: $49.987$.}
    \end{minipage}
    \hfill
    \begin{minipage}[t]{0.28\textwidth}
        \centering
        \includegraphics[width=\textwidth]{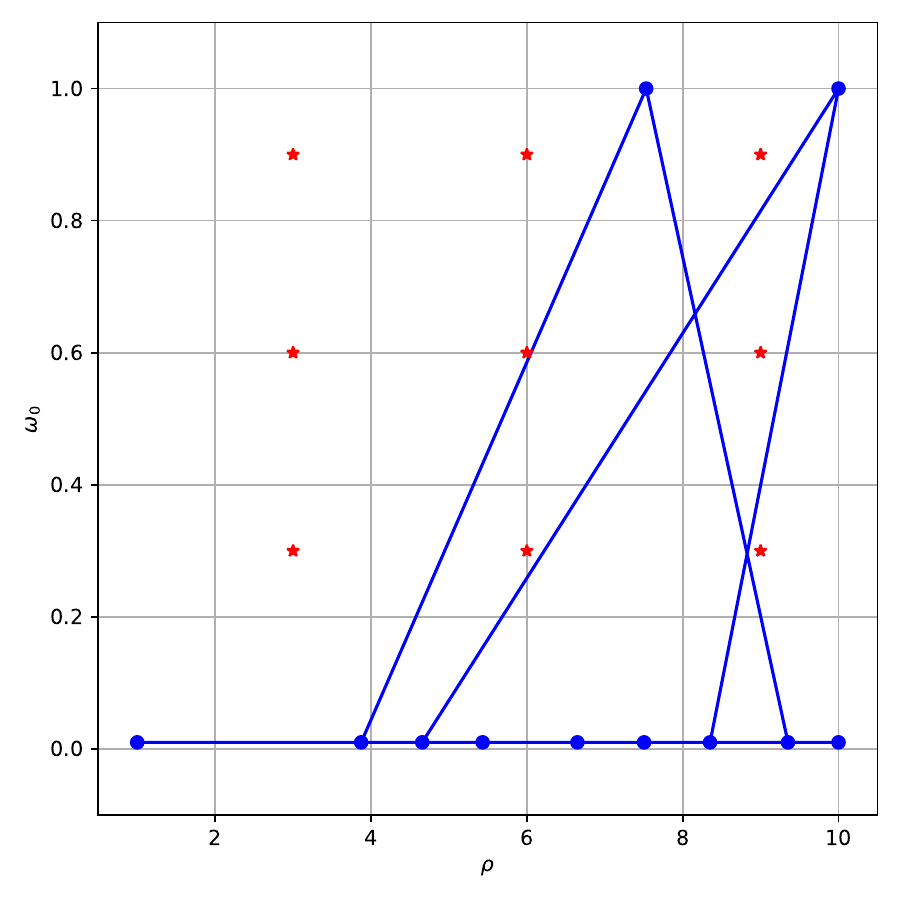}
        {\small (c) Minimum PIQUE: $49.957$.}
    \end{minipage} \vspace{0.4 cm}
    \caption{Jaszczak phantom, constant $\beta_m$, three initializations. The elements of $X_0$ are depicted in red, while in blue we show the points selected by the algorithm. The blue lines connect subsequent iterations of the scheme.}
    \label{fig:sfere_constant}
\end{figure}

\begin{figure}
    \centering
    \begin{minipage}[t]{0.28\textwidth}
        \centering
        \includegraphics[width=\textwidth]{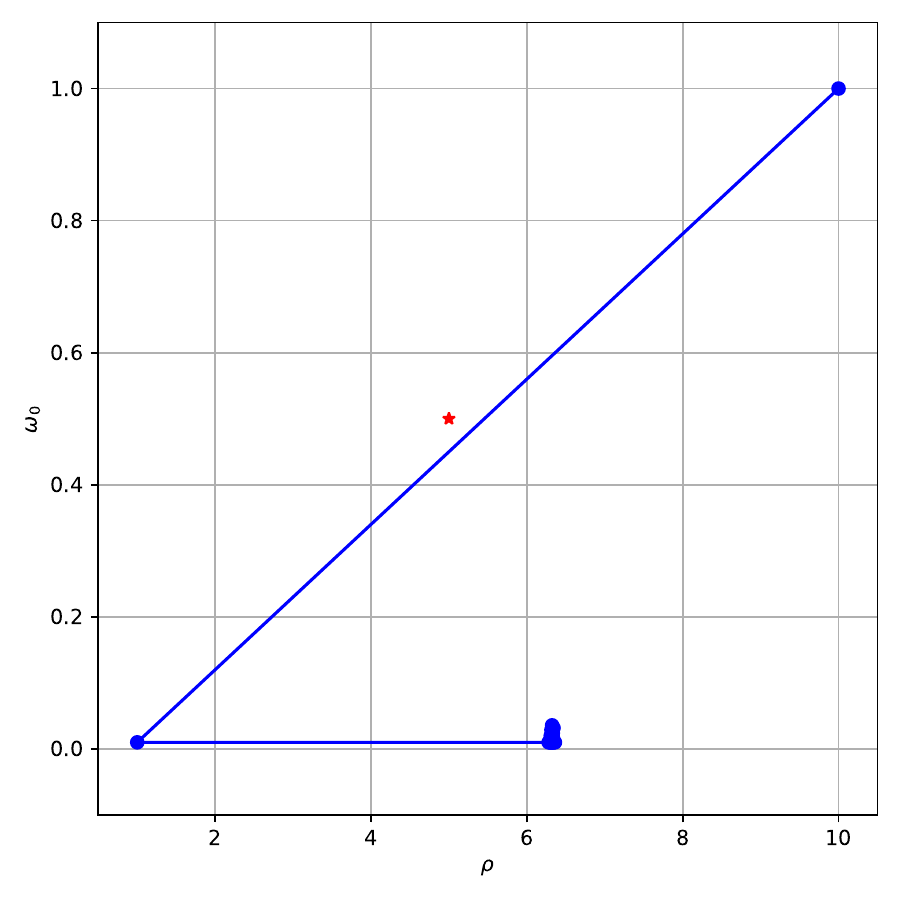}
        {\small (a) Minimum PIQUE: $50.136$.}
    \end{minipage}
    \hfill
    \begin{minipage}[t]{0.28\textwidth}
        \centering
        \includegraphics[width=\textwidth]{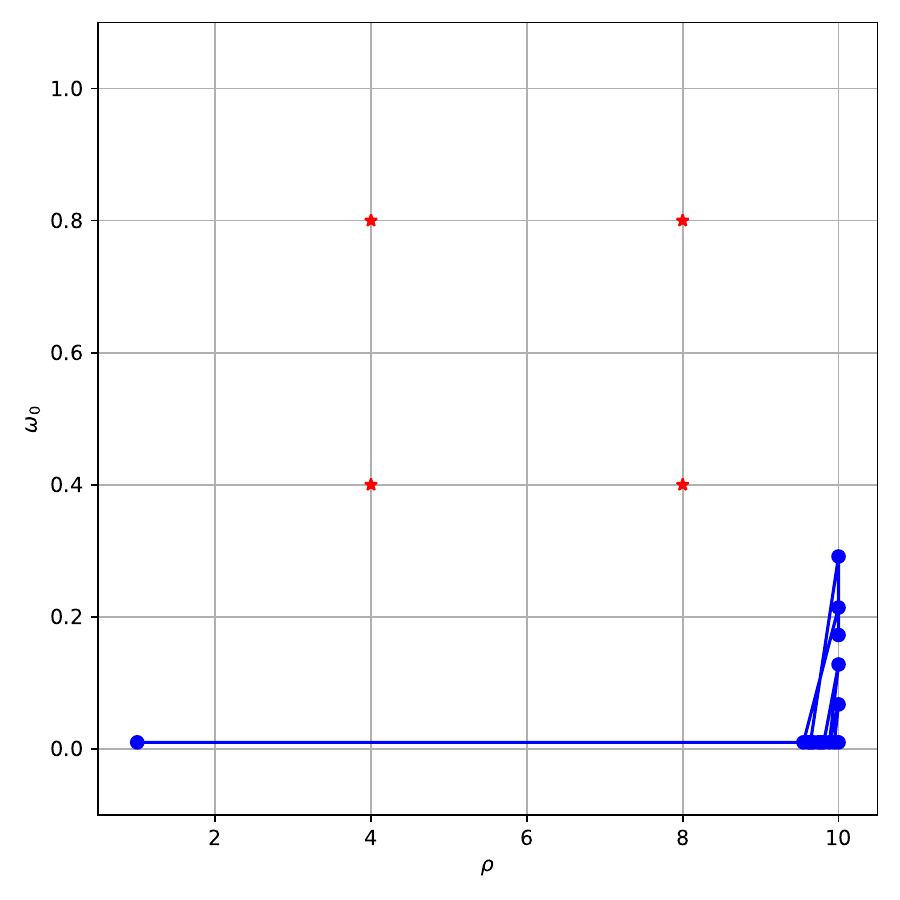}
        {\small (b) Minimum PIQUE: $49.966$.}
    \end{minipage}
    \hfill
    \begin{minipage}[t]{0.28\textwidth}
        \centering
        \includegraphics[width=\textwidth]{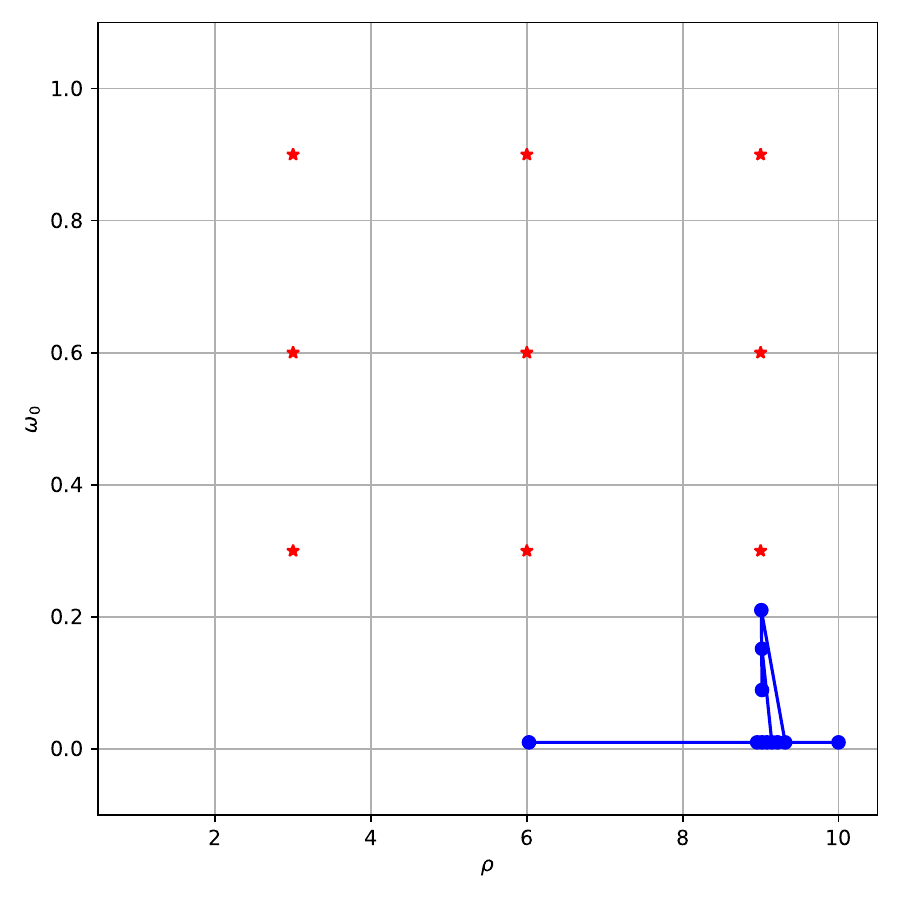}
        {\small (c) Minimum PIQUE: $49.975$.}
    \end{minipage} \vspace{0.4 cm}
    \caption{Jaszczak phantom, increasing $\beta_m$, three initializations. The elements of $X_0$ are depicted in red, while in blue we show the points selected by the algorithm. The blue lines connect subsequent iterations of the scheme.}
	\label{fig:sfere_increasing}
\end{figure}

\begin{figure}
    \centering
    \begin{minipage}[t]{0.28\textwidth}
        \centering
        \includegraphics[width=\textwidth]{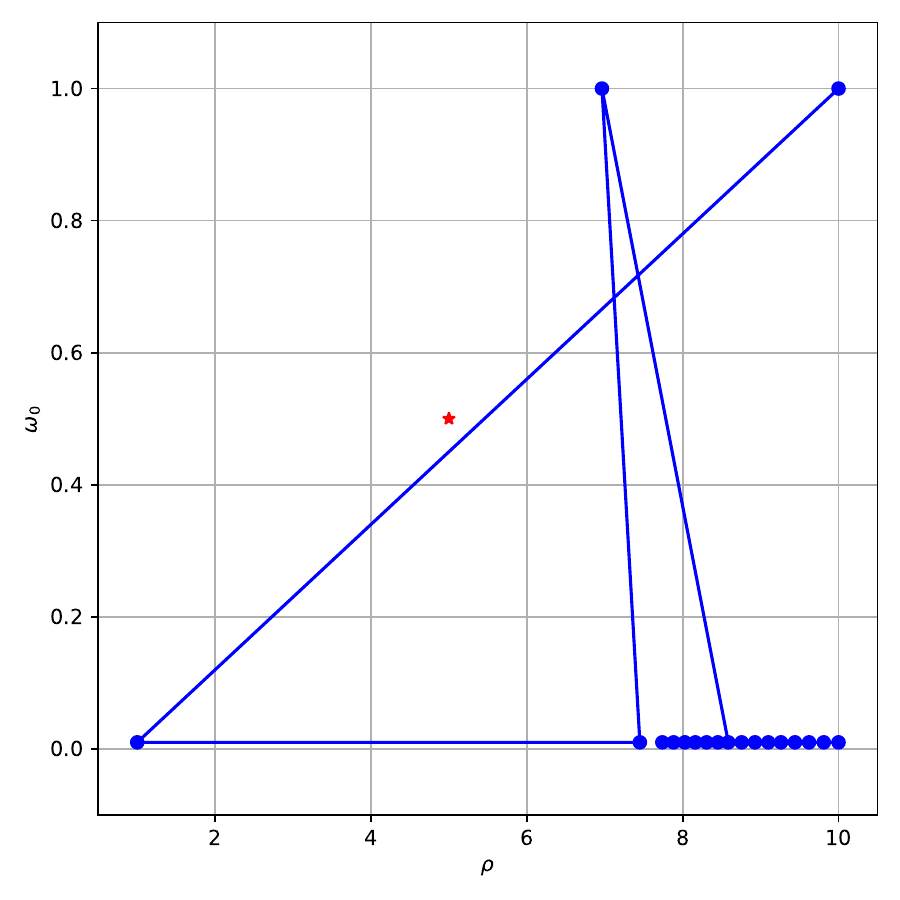}
        {\small (a) Minimum PIQUE: $49.978$.}
    \end{minipage}
    \hfill
    \begin{minipage}[t]{0.28\textwidth}
        \centering
        \includegraphics[width=\textwidth]{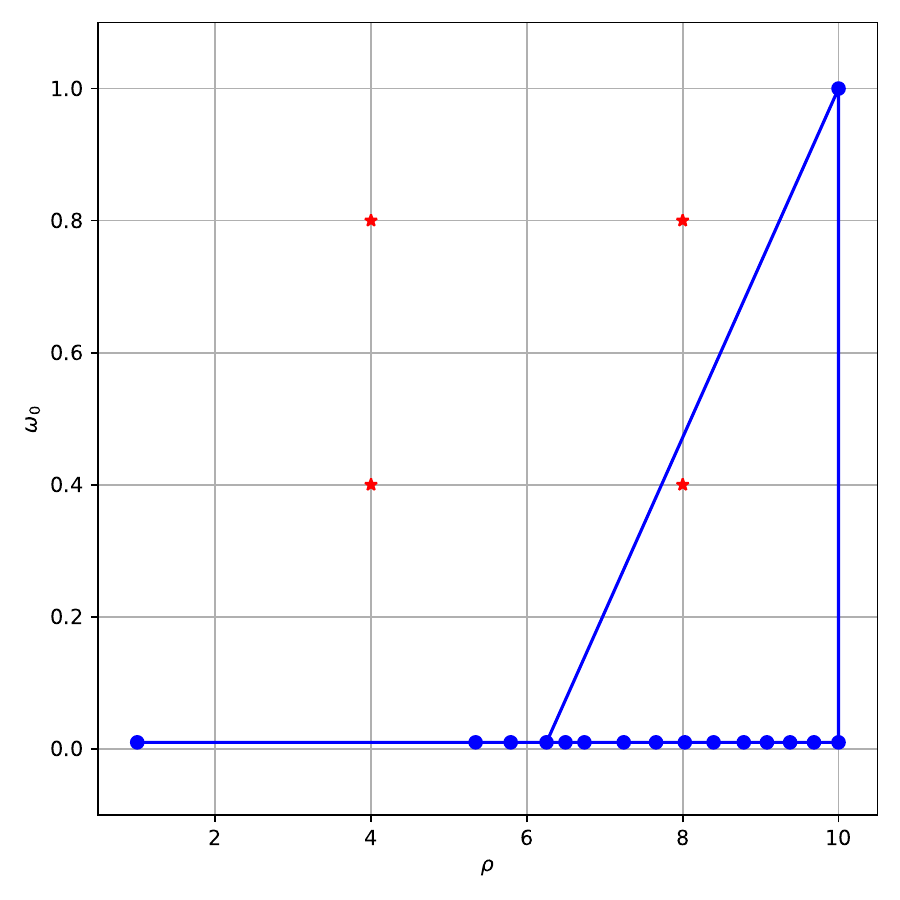}
        {\small (b) Minimum PIQUE: $49.966$.}
    \end{minipage}
    \hfill
    \begin{minipage}[t]{0.28\textwidth}
        \centering
        \includegraphics[width=\textwidth]{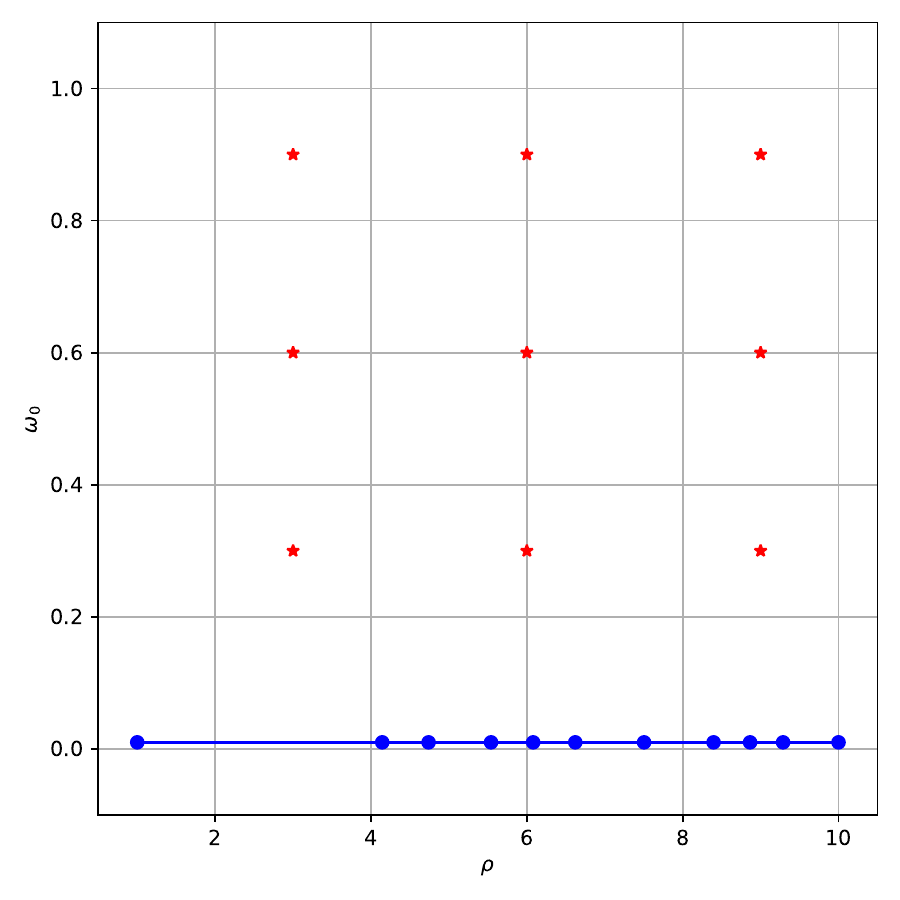}
        {\small (c) Minimum PIQUE: $49.963$.}
    \end{minipage} \vspace{0.4 cm}
    \caption{Jaszczak phantom, decreasing $\beta_m$, three initializations. The elements of $X_0$ are depicted in red, while in blue we show the points selected by the algorithm. The blue lines connect subsequent iterations of the scheme.}
	\label{fig:sfere_decreasing}
\end{figure}

Before running the BO algorithm, in Figure \ref{fig:sfere_chart} we depict the behavior of obtained PIQUE for various parameters values, and we display the FBP reconstructions corresponding to two parameter configurations in Figure \ref{fig:fbp_compare}. Then, in Figures \ref{fig:sfere_constant}, \ref{fig:sfere_increasing} and \ref{fig:sfere_decreasing}, we show the behavior of the BO greedy algorithm when considering the constant, the increasing and the decreasing weight $\beta_m$, respectively. Moreover, in these figures the three outlined initialization strategies are compared.

\subsection{Patients' data}

\begin{figure}
    \centering
    \begin{minipage}[t]{0.28\textwidth}
        \centering
        \includegraphics[width=\textwidth]{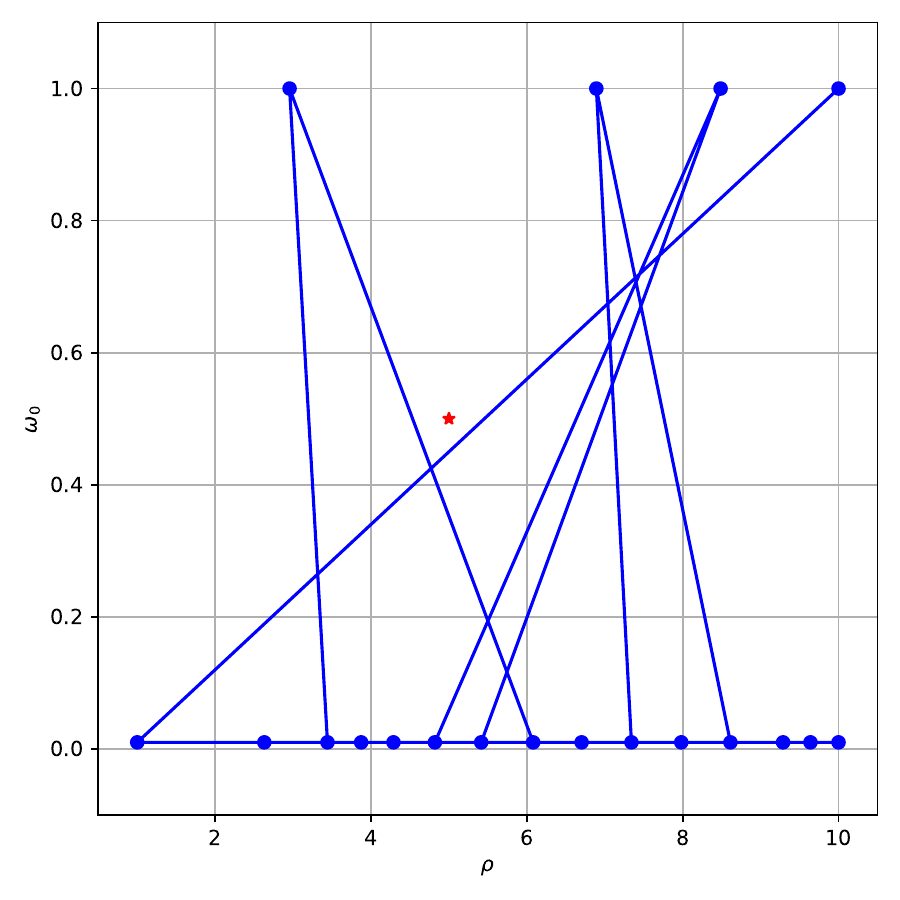}
        {\small (a) Minimum PIQUE: $35.605$.}
    \end{minipage}
    \hfill
    \begin{minipage}[t]{0.28\textwidth}
        \centering
        \includegraphics[width=\textwidth]{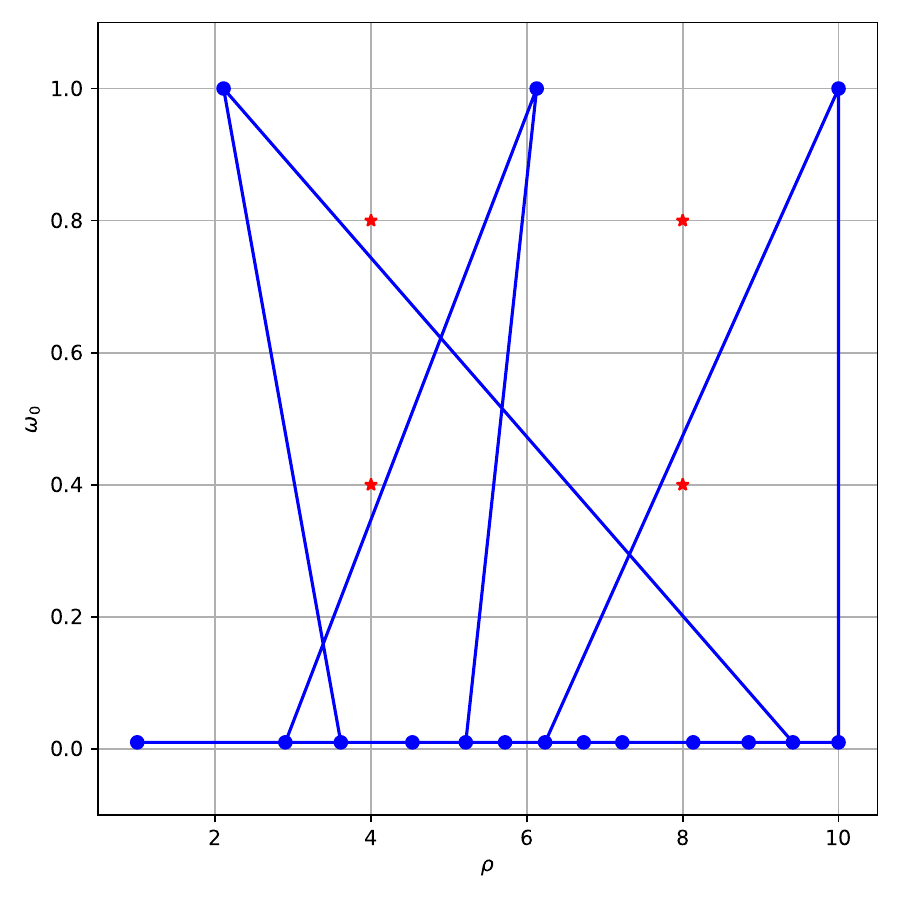}
        {\small (b) Minimum PIQUE: $35.609$.}
    \end{minipage}
    \hfill
    \begin{minipage}[t]{0.28\textwidth}
        \centering
        \includegraphics[width=\textwidth]{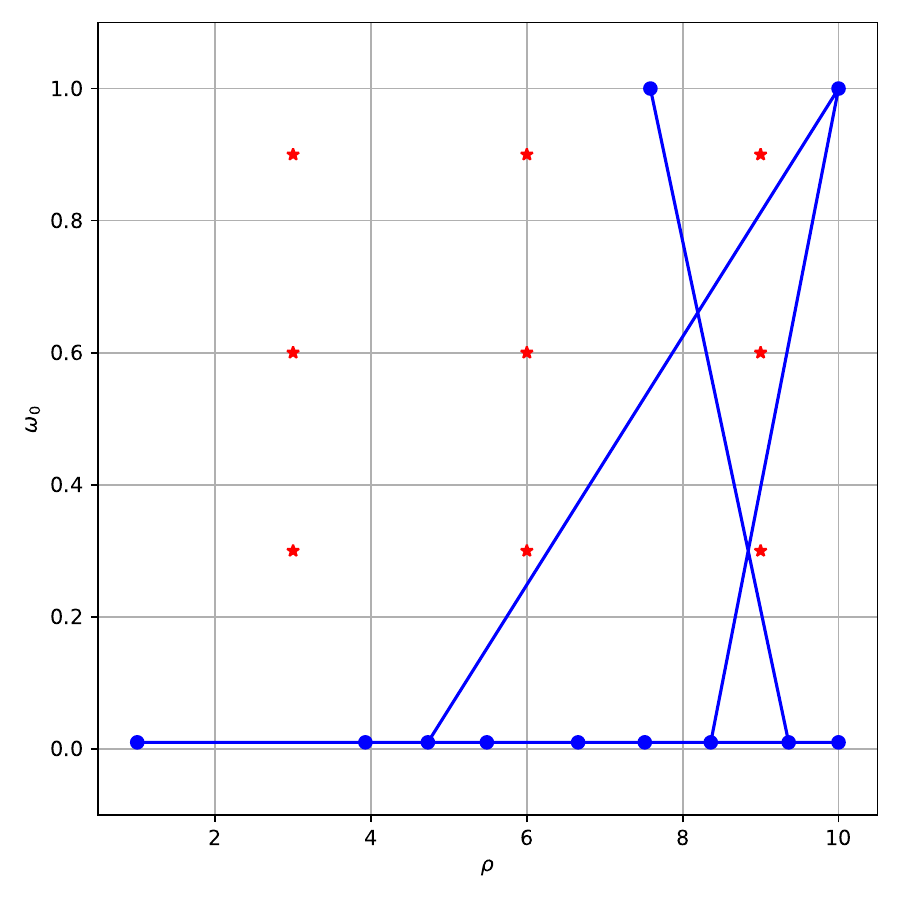}
        {\small (c) Minimum PIQUE: $35.629$.}
    \end{minipage} \vspace{0.4 cm}
    \caption{Lung perfusion scintigraphy, constant $\beta_m$, three initializations. The elements of $X_0$ are depicted in red, while in blue we show the points selected by the algorithm. The blue lines connect subsequent iterations of the scheme.}
	\label{fig:polmone_constant}
\end{figure}

\begin{figure}
    \centering
    \begin{minipage}[t]{0.28\textwidth}
        \centering
        \includegraphics[width=\textwidth]{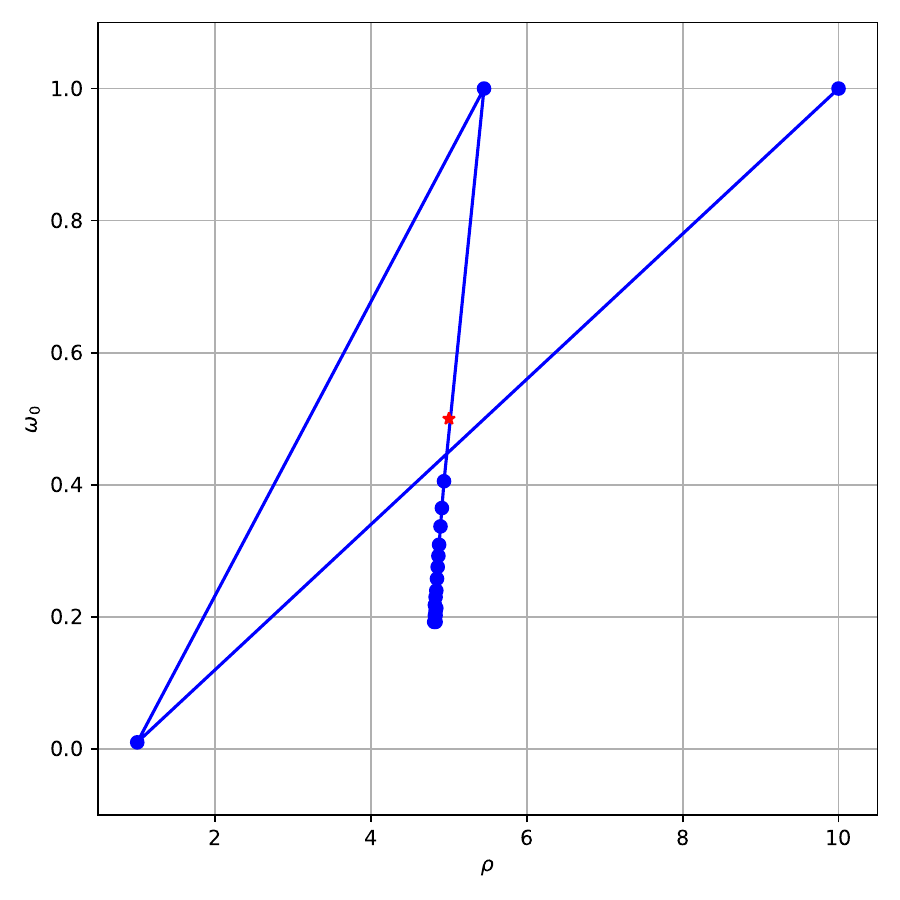}
        {\small (a) Minimum PIQUE: $36.106$.}
    \end{minipage}
    \hfill
    \begin{minipage}[t]{0.28\textwidth}
        \centering
        \includegraphics[width=\textwidth]{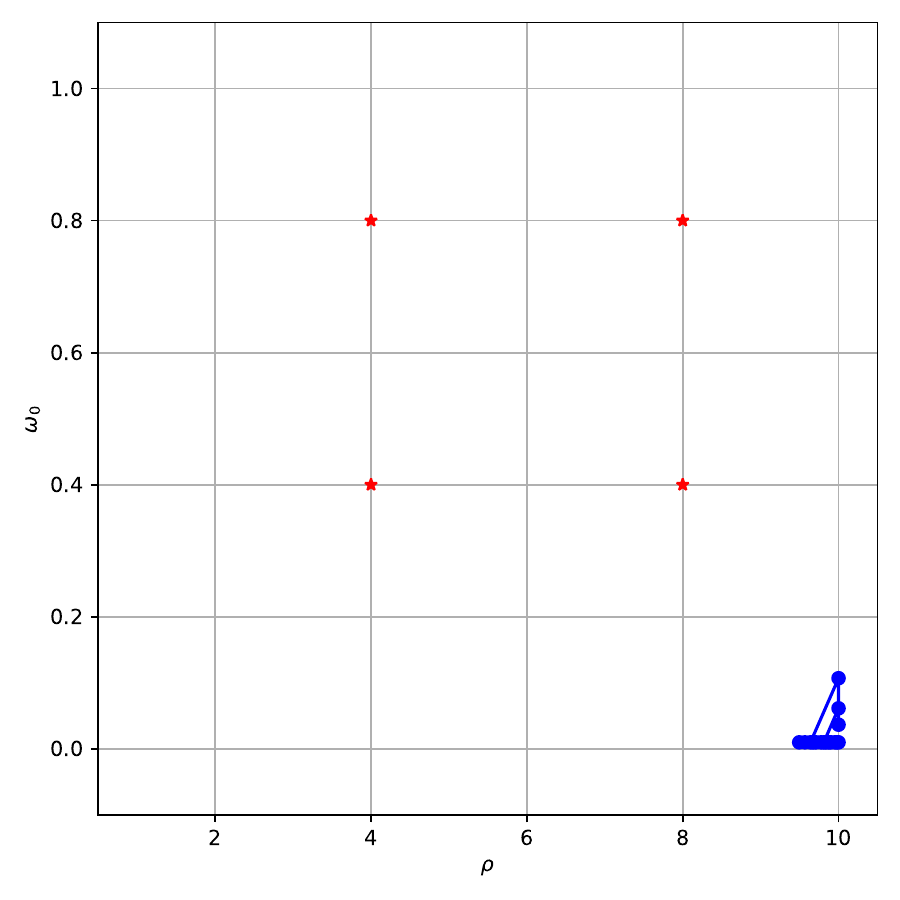}
        {\small (b) Minimum PIQUE: $35.618$.}
    \end{minipage}
    \hfill
    \begin{minipage}[t]{0.28\textwidth}
        \centering
        \includegraphics[width=\textwidth]{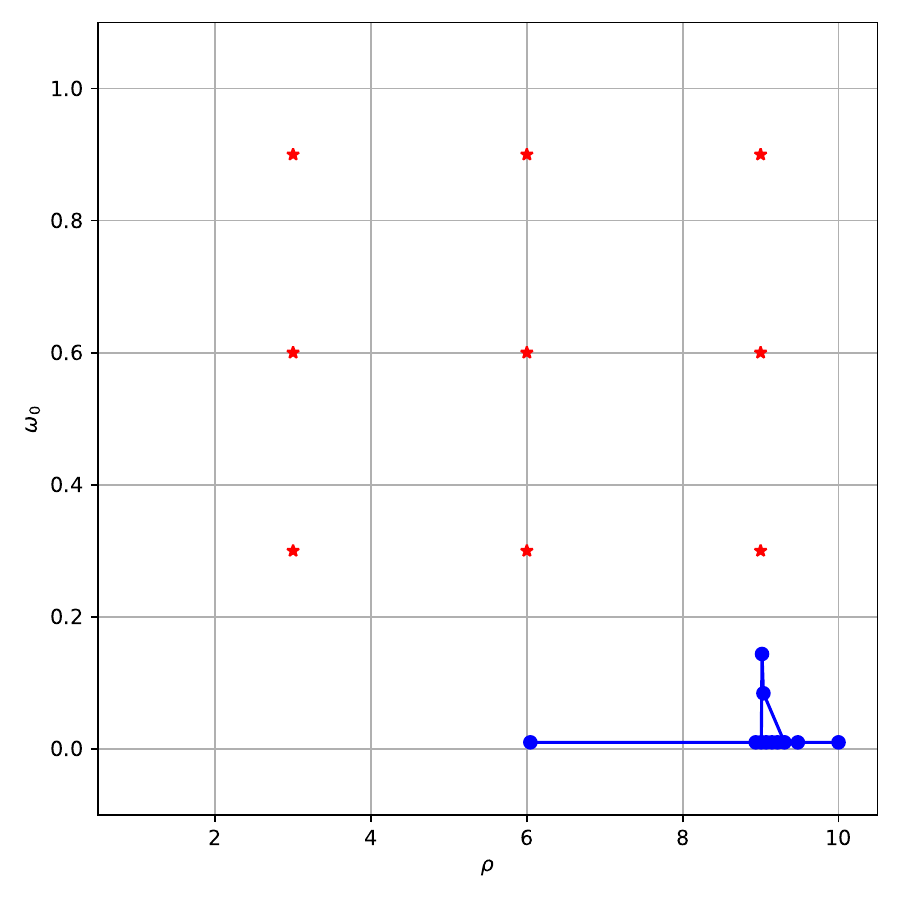}
        {\small (c) Minimum PIQUE: $35.601$.}
    \end{minipage} \vspace{0.4 cm}
    \caption{Lung perfusion scintigraphy, increasing $\beta_m$, three initializations. The elements of $X_0$ are depicted in red, while in blue we show the points selected by the algorithm. The blue lines connect subsequent iterations of the scheme.}
	\label{fig:polmone_increasing}
\end{figure}

\begin{figure}
    \centering
    \begin{minipage}[t]{0.28\textwidth}
        \centering
        \includegraphics[width=\textwidth]{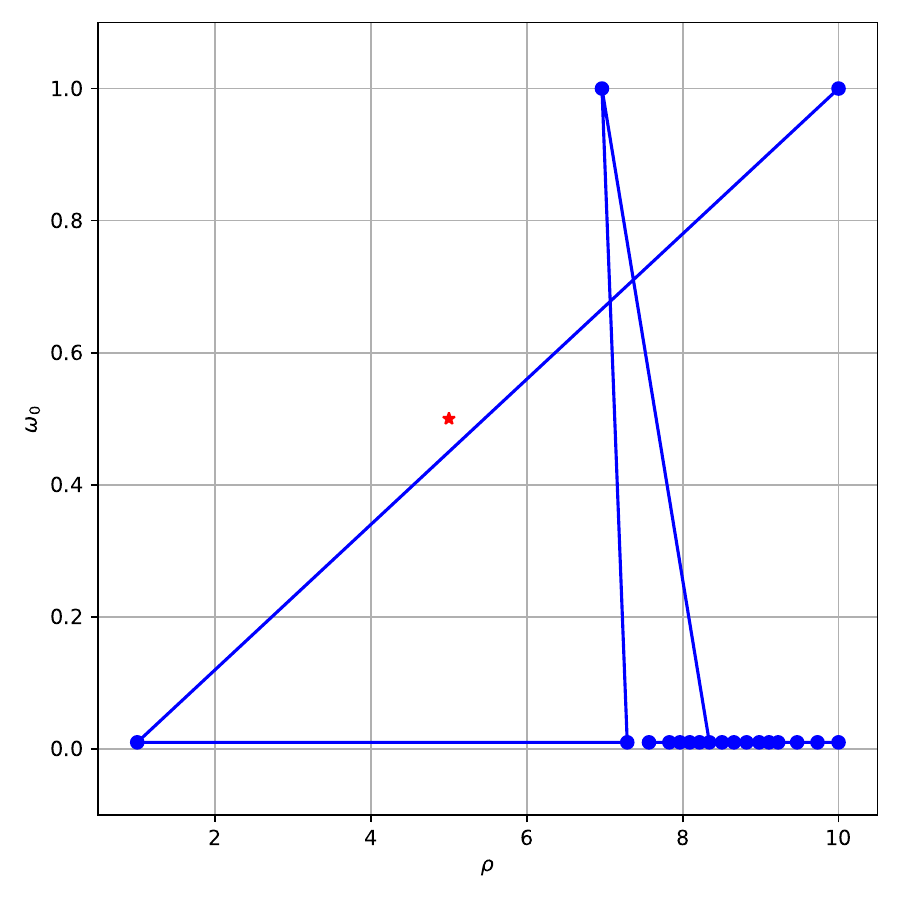}
        {\small (a) Minimum PIQUE: $35.589$.}
    \end{minipage}
    \hfill
    \begin{minipage}[t]{0.28\textwidth}
        \centering
        \includegraphics[width=\textwidth]{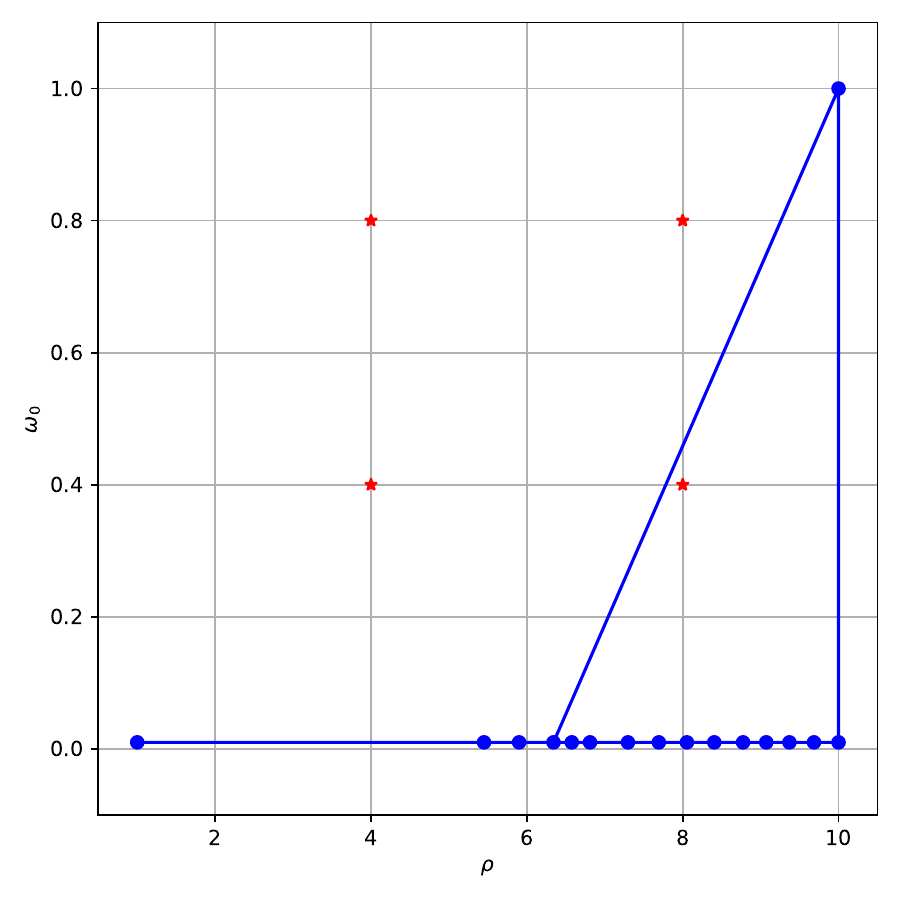}
        {\small (b) Minimum PIQUE: $35.604$.}
    \end{minipage}
    \hfill
    \begin{minipage}[t]{0.28\textwidth}
        \centering
        \includegraphics[width=\textwidth]{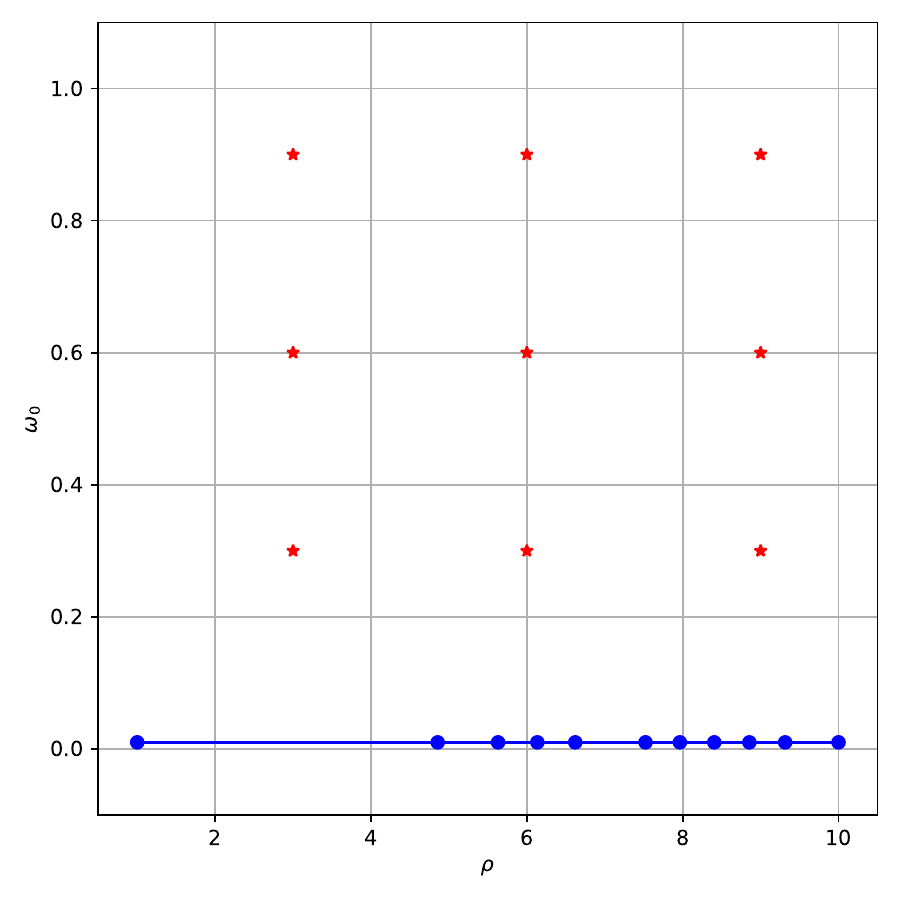}
        {\small (c) Minimum PIQUE: $35.604$.}
    \end{minipage} \vspace{0.4 cm}
    \caption{Lung perfusion scintigraphy, decreasing $\beta_m$, three initializations. The elements of $X_0$ are depicted in red, while in blue we show the points selected by the algorithm. The blue lines connect subsequent iterations of the scheme.}
	\label{fig:polmone_decreasing}
\end{figure}

\begin{figure}
    \centering
    \begin{minipage}[t]{0.28\textwidth}
        \centering
        \includegraphics[width=\textwidth]{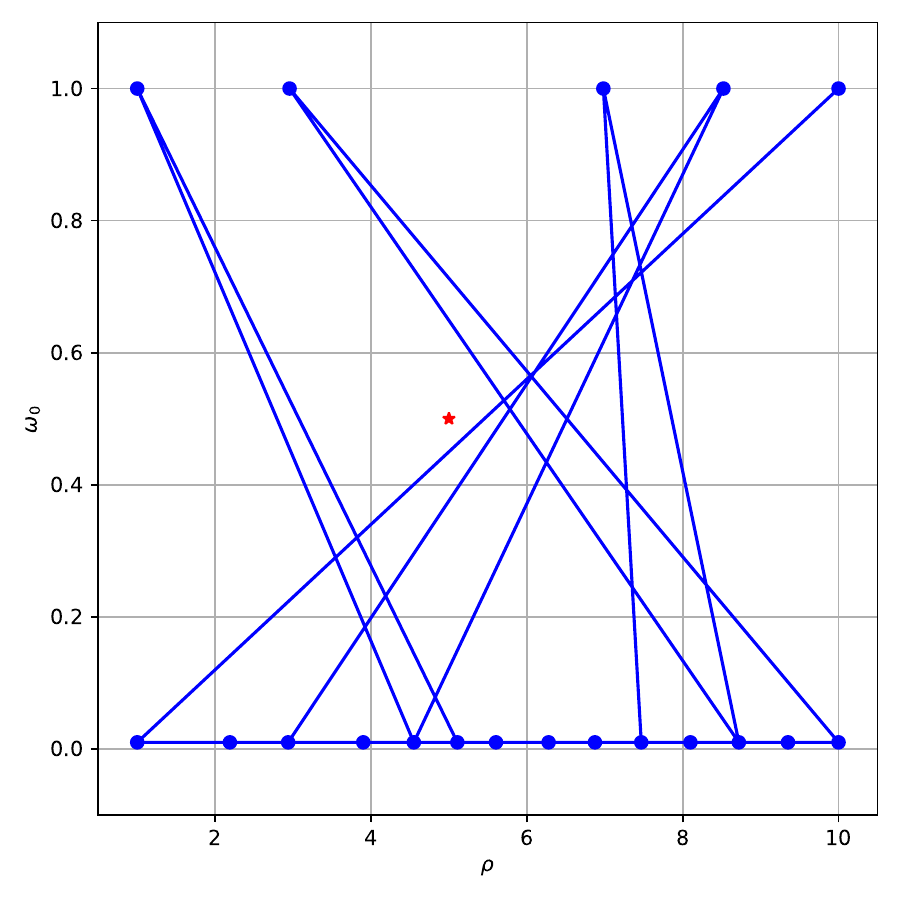}
        {\small (a) Minimum PIQUE: $42.614$.}
    \end{minipage}
    \hfill
    \begin{minipage}[t]{0.28\textwidth}
        \centering
        \includegraphics[width=\textwidth]{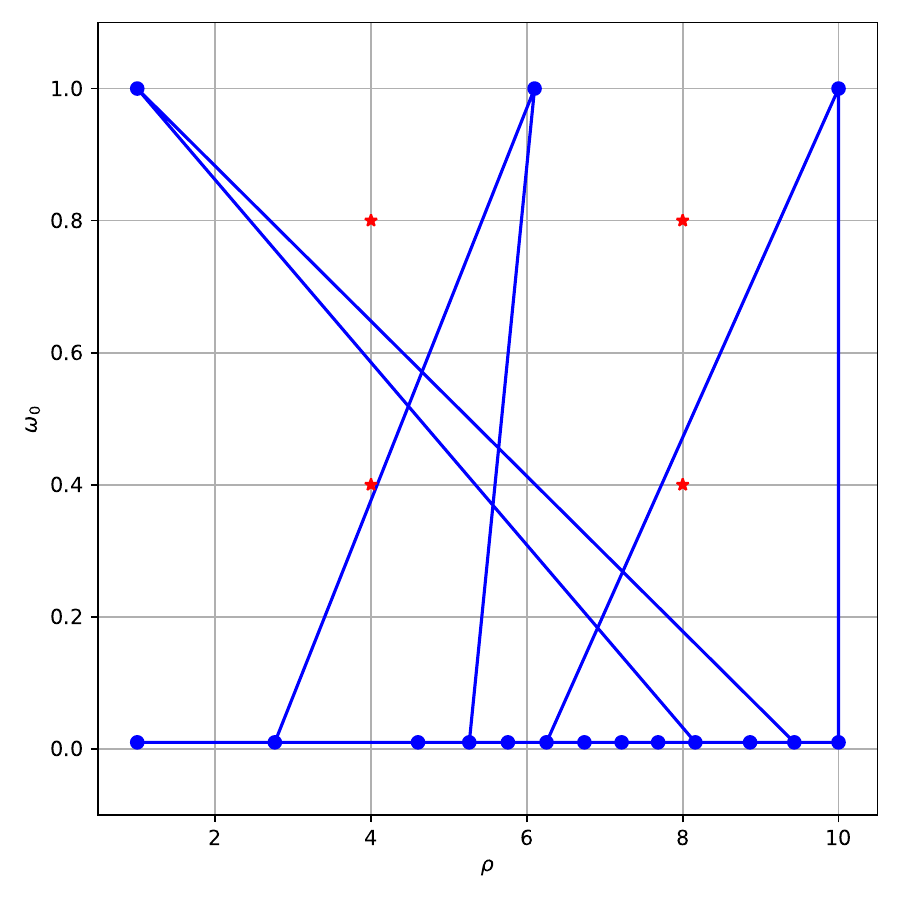}
        {\small (b) Minimum PIQUE: $42.614$.}
    \end{minipage}
    \hfill
    \begin{minipage}[t]{0.28\textwidth}
        \centering
        \includegraphics[width=\textwidth]{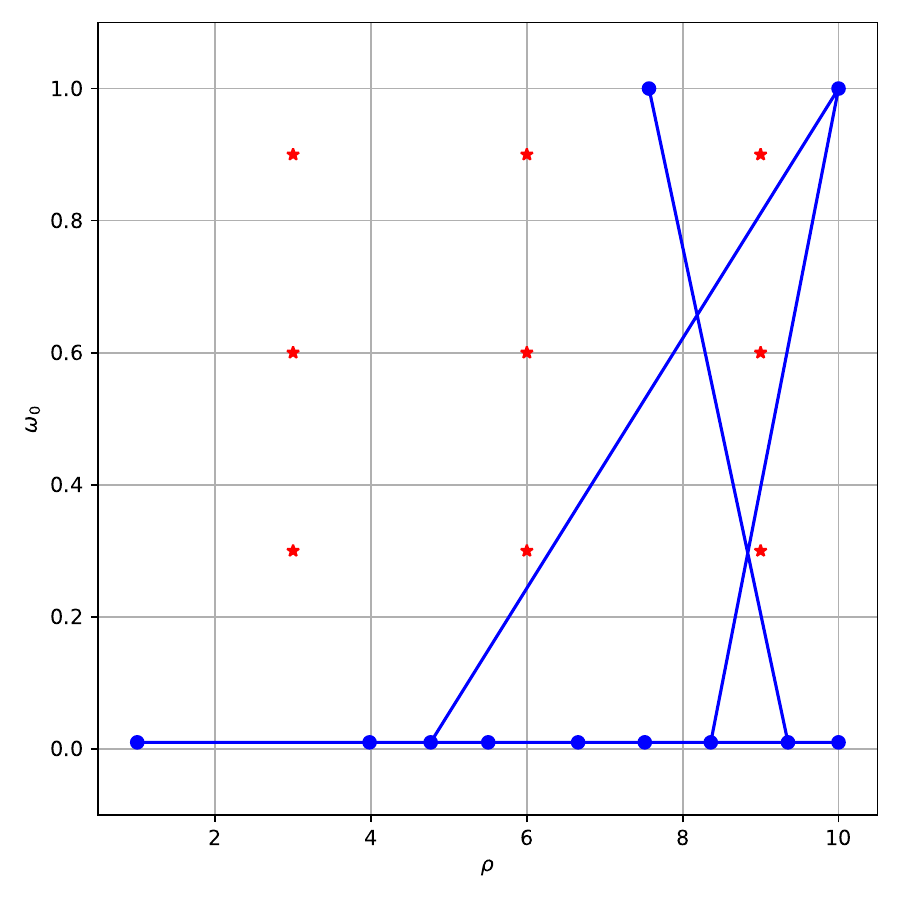}
        {\small (c) Minimum PIQUE: $42.614$.}
    \end{minipage} \vspace{0.4 cm}
    \caption{Thyroid scintigraphy, constant $\beta_m$, three initializations. The elements of $X_0$ are depicted in red, while in blue we show the points selected by the algorithm. The blue lines connect subsequent iterations of the scheme.}
	\label{fig:tiroide_constant}
\end{figure}

\begin{figure}
    \centering
    \begin{minipage}[t]{0.28\textwidth}
        \centering
        \includegraphics[width=\textwidth]{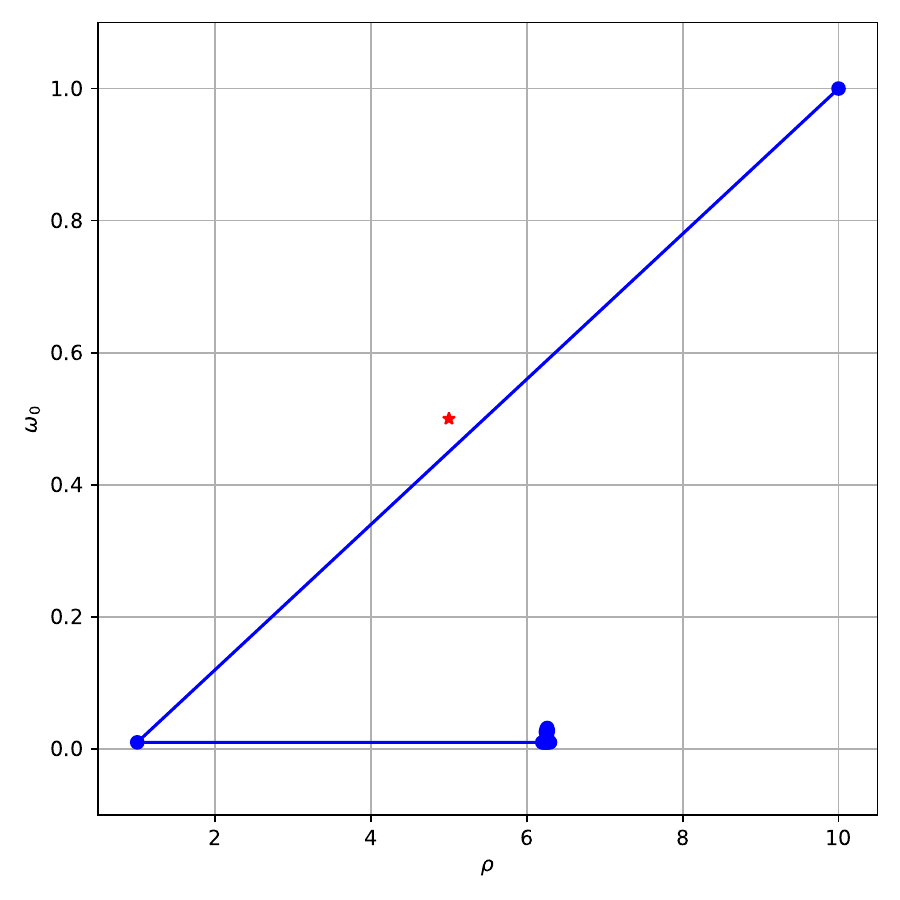}
        {\small (a) Minimum PIQUE: $42.752$.}
    \end{minipage}
    \hfill
    \begin{minipage}[t]{0.28\textwidth}
        \centering
        \includegraphics[width=\textwidth]{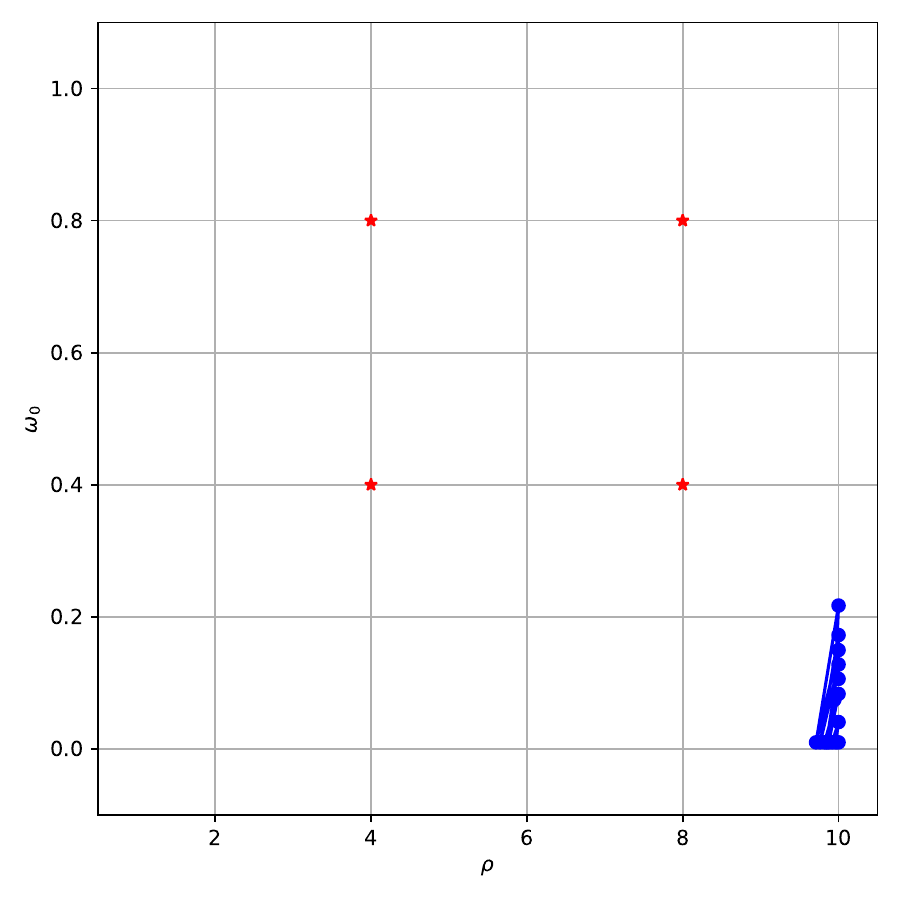}
        {\small (b) Minimum PIQUE: $42.604$.}
    \end{minipage}
    \hfill
    \begin{minipage}[t]{0.28\textwidth}
        \centering
        \includegraphics[width=\textwidth]{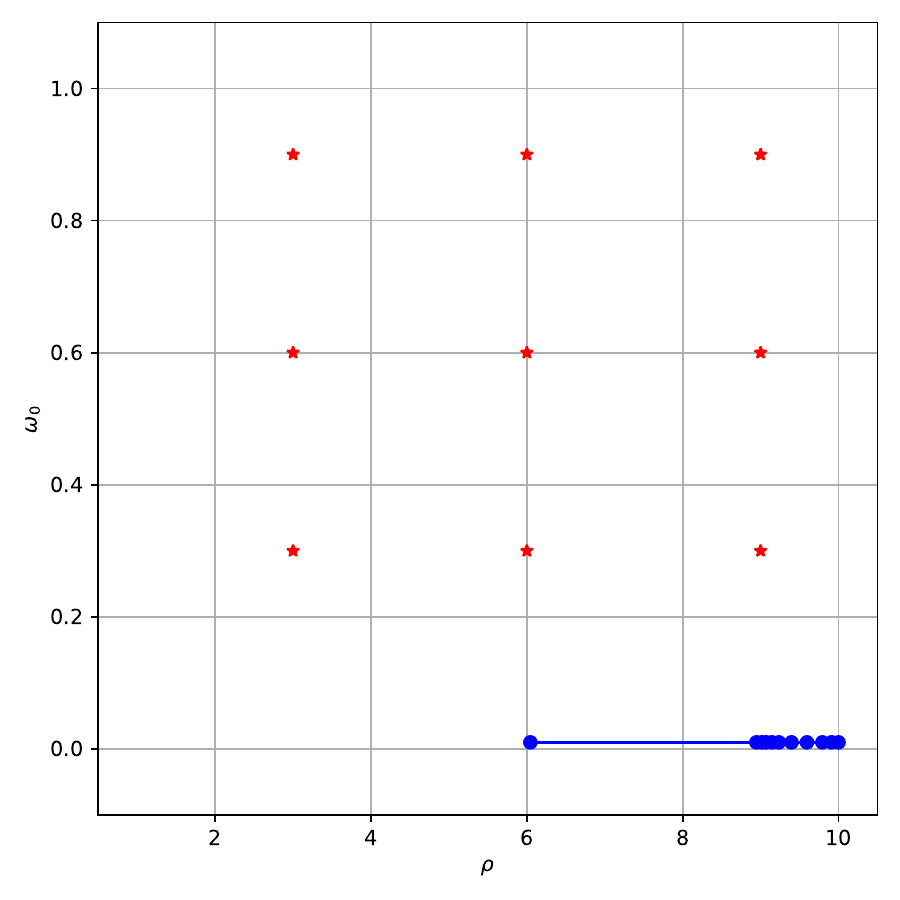}
        {\small (c) Minimum PIQUE: $42.614$.}
    \end{minipage} \vspace{0.4 cm}
    \caption{Thyroid scintigraphy, increasing $\beta_m$, three initializations. The elements of $X_0$ are depicted in red, while in blue we show the points selected by the algorithm. The blue lines connect subsequent iterations of the scheme.}
	\label{fig:tiroide_increasing}
\end{figure}

\begin{figure}
    \centering
    \begin{minipage}[t]{0.28\textwidth}
        \centering
        \includegraphics[width=\textwidth]{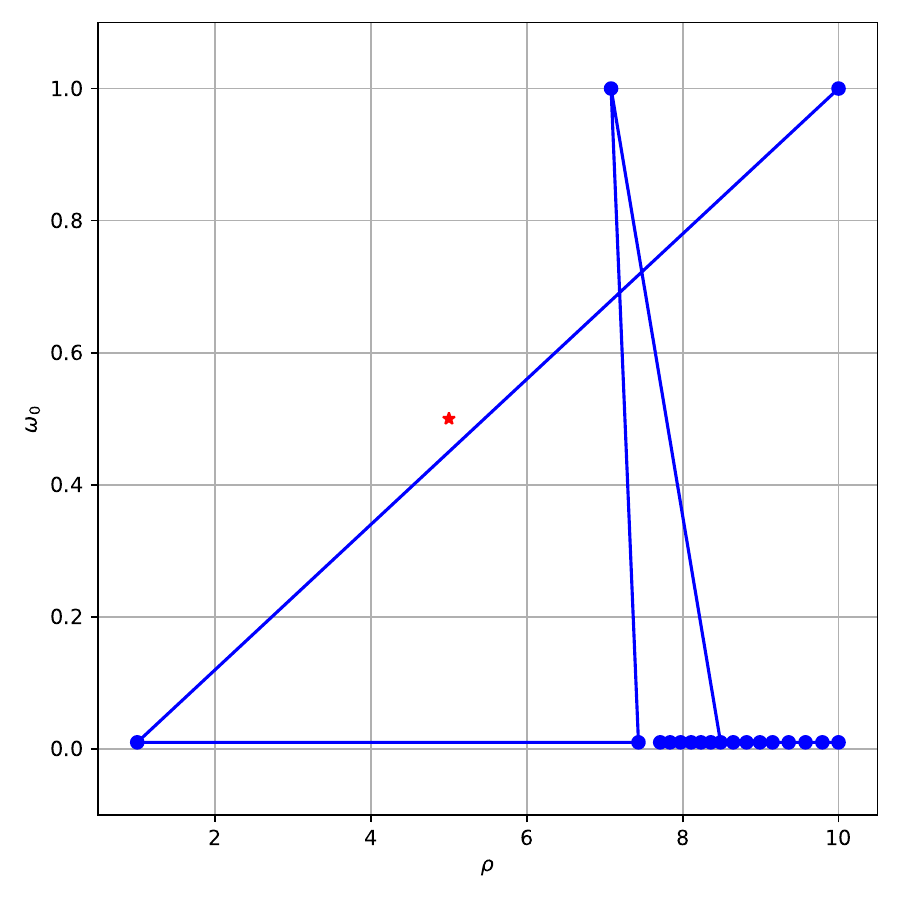}
        {\small (a) Minimum PIQUE: $42.614$.}
    \end{minipage}
    \hfill
    \begin{minipage}[t]{0.28\textwidth}
        \centering
        \includegraphics[width=\textwidth]{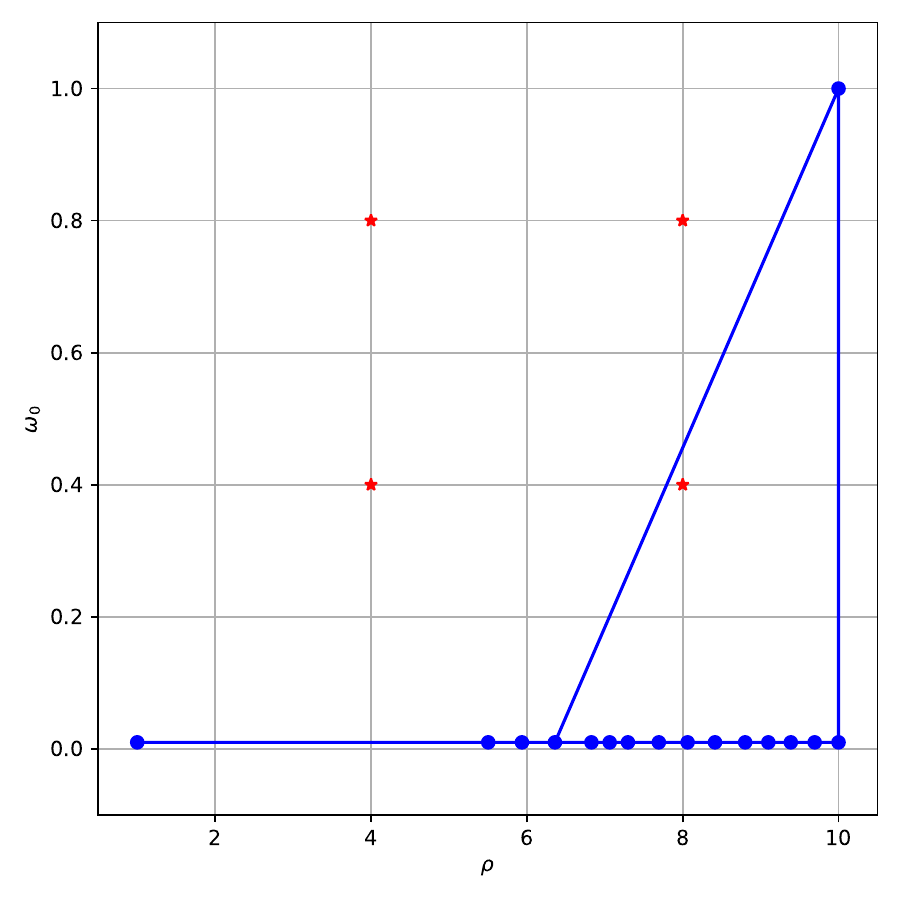}
        {\small (b) Minimum PIQUE: $42.614$.}
    \end{minipage}
    \hfill
    \begin{minipage}[t]{0.28\textwidth}
        \centering
        \includegraphics[width=\textwidth]{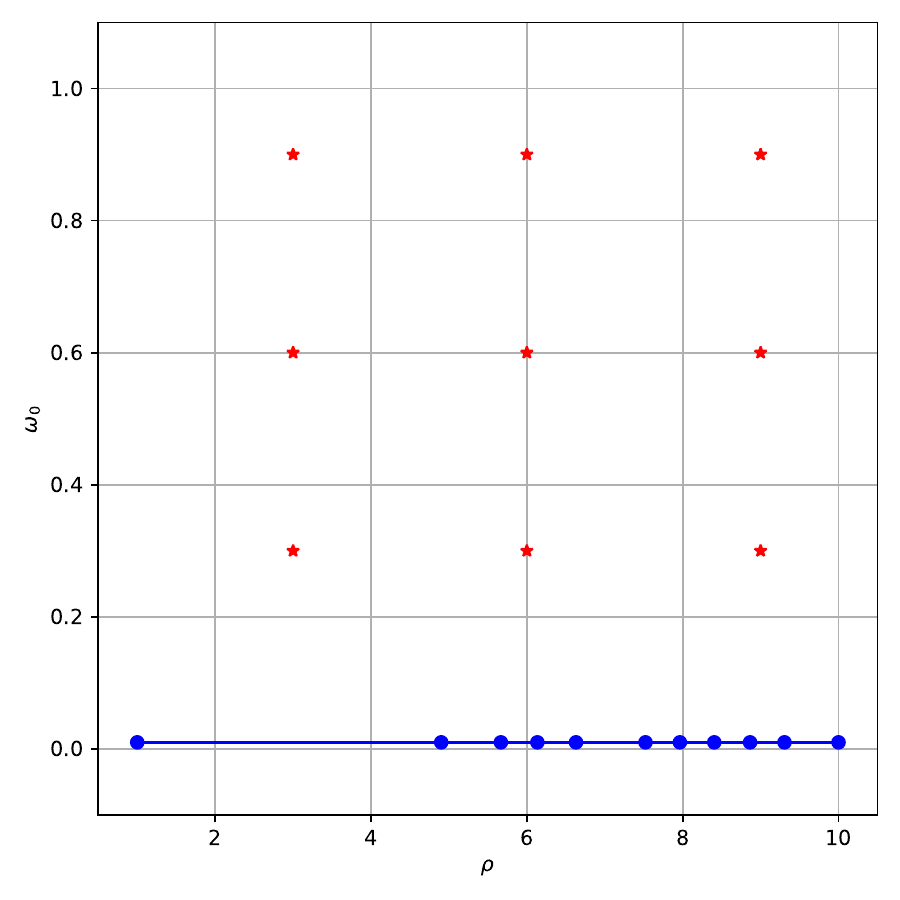}
        {\small (c) Minimum PIQUE: $42.614$.}
    \end{minipage} \vspace{0.4 cm}
    \caption{Thyroid scintigraphy, decreasing $\beta_m$, three initializations. The elements of $X_0$ are depicted in red, while in blue we show the points selected by the algorithm. The blue lines connect subsequent iterations of the scheme.}
	\label{fig:tiroide_decreasing}
\end{figure}

In this subsection, we experiment with two acquisitions obtained in concrete clinical settings.
\begin{itemize}
    \item 
    Lung perfusion scintigraphy: 99mTc-MAA (approximately 120 MBq). Acquired approximately 5 minutes after tracer injection, which was administered slowly while the patient did deep inspiration and expiration.
    \item 
    Thyroid scintigraphy: 99mTc (approximately 150 MBq). Acquired about 15 minutes after tracer injection.
\end{itemize}
 Similarly to the previous subsection, the results for lung perfusion scintigraphy are shown in Figures \ref{fig:polmone_constant}, \ref{fig:polmone_increasing} and \ref{fig:polmone_decreasing}, while those for the thyroid are shown in Figures \ref{fig:tiroide_constant}, \ref{fig:tiroide_increasing} and \ref{fig:tiroide_decreasing}.

\section{Discussion and conclusions}

First, we observe that Figure \ref{fig:fbp_compare} is an indication of the fact that the no-reference PIQUE metric is meaningful in SPECT image reconstruction tasks, as already observed for PET in previous works. Indeed, the usage of optimal parameters configuration leads to a better-defined image that contains a smaller amount of background noise.

The optimization task is well addressed by kernel-based BO. On the one hand, as presented in Figure \ref{fig:sfere_chart}, there exists a \textit{plateau} of (quasi-)optimal configurations that facilitate a fine-tuning of the filter. On the other hand, we observe that some \textit{discontinuous} behavior is shown in the figure under some circumstances; therefore, using a gradient-free approach is a reasonable choice. We also remark that the phantom's framework and the patients' acquisition share very similar good configurations of parameters in terms of the chosen PIQUE metric, which is crucial in order to extend a fine-tuned setting to other image acquisition tasks. Unlike the other approaches, the use of increasing weights as in Figures \ref{fig:sfere_increasing}, \ref{fig:polmone_increasing} and \ref{fig:tiroide_increasing} limits the space exploration and leads to a fairly well-outlined optimization path that, however, results in being \textit{slower} than the other schemes.

Therefore, the results of this paper can stimulate future development in the usage of the PIQUE metric and kernel-based greedy optimization strategies in SPECT imaging, including iterative approaches. From a theoretical point of view, an interesting research line to ensure a robust framework could be the analysis of convergence rates for PIQUE similar to what was carried out for SSIM in \cite{Marchetti22}.

\section*{Acknowledgments}
This research has been accomplished within Rete ITaliana di Approssimazione (RITA), with the support of GNCS-IN$\delta$AM.


\begin{thebibliography}{10}


  
 	\bibitem{Bailey05}
 \textsc{D.L. Bailey, D.W. Townsend, P.E. Valk, M.N. Maisey},
 \emph{Positron Emission Tomography: Basic Sciences}, 
 Springer, 2005.
 
  	\bibitem{Bull11}
 \textsc{A. D. Bull}, \emph{Convergence rates of efficient global optimization algorithms}, JMLR \textbf{12} (2011), pp. 2879--2904.
 
 	\bibitem{Cavoretto18}
 \textsc{R. Cavoretto, A. De Rossi, E. Perracchione}, \emph{Optimal Selection of Local Approximants in RBF-PU Interpolation}, J Sci Comput \textbf{74} (2018), pp. 1--22.
 
 \bibitem{Chen22}
 \textsc{Y.-T. Chen, C. Li, L.-Q. Yao, Y. Cao}, \emph{A Hybrid RBF Collocation Method and Its Application in the Elastostatic Symmetric Problems}, Symmetry \textbf{14}(7) (2022), 1476.
 
	\bibitem{Cherry12}
\textsc{S.R. Cherry, J.A. Sorenson, M.E. Phelps},
\emph{Physics in Nuclear Medicine (4th ed.)}, 
Elsevier, 2012.
 
 	\bibitem{DeMarchi05} 
 \textsc{S. De Marchi, R. Schaback, H. Wendland}, 
 \emph{Near-optimal data-independent point locations for radial basis function interpolation}, Adv. Comp. Math. \textbf{23}(3) (2005), pp. 317--330.


	\bibitem{Fasshauer07}
\textsc{G.E. Fasshauer},
\emph{Meshfree Approximations Methods with \textsc{Matlab}}, 
World Scientific, Singapore, 2007.

\bibitem{Fasshauer15}
\textsc{G.E. Fasshauer, M.J. McCourt}, 
\emph{Kernel-based Approximation Methods Using \textsc{Matlab}}, 
World Scientific, Singapore, 2015.

	\bibitem{Fornberg07} 
\textsc{B. Fornberg, J. Zuev}, \emph{The Runge phenomenon and spatially variable shape parameters in RBF interpolation}, Comput. Math. Appl. \textbf{54}(3) (2007), pp. 379--398.

\bibitem{shigeaki2024}
\textsc{S. Higashiyama, T. Yamanaga, Y. Katayama, A. Yoshida, N. Inoue, T. Ichida, Y. Miki, J. Kawabe}, 
\emph{Investigation of the effectiveness of no reference metric in image evaluation in nuclear medicine}, PLOS ONE 2024.



	\bibitem{Hutton19} 
\textsc{B.F. Hutton, K. Erlandsson}, \emph{Advances in clinical molecular imaging with SPECT/CT and PET/CT}, Clinical Physiology and Functional Imaging \textbf{39}(2) (2019), pp. 81--95.

	\bibitem{Hutton11} 
\textsc{B.F. Hutton, J. Nuyts, F.J. Beekman}, 
\emph{Review of image reconstruction methods for PET and SPECT}, The European Journal of Nuclear Medicine and Molecular Imaging \textbf{38}(1) (2011), pp. 167--188.	



\bibitem{Jaszczak80} 
\textsc{R.J. Jaszczak et al.}, 
\emph{Physical Performance Characteristics of a Whole-Body Tomographic Scanner}, J. Nucl. Med. \textbf{21}(4) (1980), pp. 344–-350.	
  .
  
\bibitem{Kak01}
\textsc{A.C. Kak, M. Slaney}, 
\emph{Principles of Computerized Tomographic Imaging}, 
SIAM, 2001.

\bibitem{kanagawa2018}
\textsc{M. Kanagawa, P. Hennig, D. Sejdinovic, B.K. Sriperumbudur},
\emph{ Gaussian processes and kernel methods: A review on connections and equivalences}, arXiv:1807.02582, 2018.

  
  
  \bibitem{LeRouzic21} 
  \textsc{G. Le Rouzic, R. Zananiri}, 
  \emph{First performance measurements of a new multi-detector CZT-Based SPECT/CT system: GE StarGuide}, J. Nucl. Med. \textbf{62}(1) (2021), pp. 1125–-1125.	
  .
  \bibitem{Lyu19} 
  \textsc{Y. Lyu, Y. Yuan, I.W. Tsang}, 
  \emph{Efficient batch black-box optimization with deterministic regret bounds}, arXiv:1905.10041, 2019.
  
  
  \bibitem{Marchetti24}
  \textsc{F. Marchetti}, \emph{A fast surrogate cross validation algorithm for meshfree RBF collocation approaches}, Appl. Math. Comput. \textbf{481} (2024), 128943.
  
  \bibitem{Marchetti22}
  \textsc{F. Marchetti, G. Santin}, \emph{Convergence results in image interpolation with the continuous SSIM}, SIAM J. Imag. Sci. \textbf{15}(4) (2022), pp. 1977--1999.
  

\bibitem{mittal2012}
\textsc{A. Mittal, A. Moorthy, A. Bovik}
\emph{No-reference image quality assessment in the spatial domain}, Transactions on Image Processing, IEEE \textbf{21}(12) (2012), pp. 4695--4708. 

	
	
	\bibitem{Narcowich05} 
	\textsc{F.J. Narcowich, J.D. Ward, H. Wendland}, 
	\emph{Sobolev bounds on functions with scattered zeros, with applications to radial basis function surface fitting}, Math. Comp. \textbf{74} (2005), pp. 743--763.
	
	\bibitem{Natterer86}
	\textsc{F. Natterer}, 
	\emph{The Mathematics of Computerized Tomography}, 
	SIAM, 1986.
	
	
	
	
	\bibitem{Rahmim08} 
	\textsc{A. Rahmim, H. Zaidi}, 
	\emph{PET versus SPECT: Strengths, limitations, and challenges}, Nuclear Medicine Communications \textbf{29}(3) (2015), pp. 193--207.
	
	
	\bibitem{Reader07} 
	\textsc{A.J. Reader, H. Zaidi}, 
	\emph{Advances in PET image reconstruction}, PET Clinics \textbf{2}(2) (2007), pp. 173--190.	
	
	\bibitem{Rieger09} 
	\textsc{C. Rieger, B. Zwicknagl}, 
	\emph{Deterministic Error Analysis of Support Vector Regression and Related Regularized Kernel Methods}, J. Mach. Learn. Res. \textbf{10} (2009), pp. 2115--2132.
	
	
	
	\bibitem{Santin16} 
	\textsc{G. Santin, B. Haasdonk}, 
	\emph{Convergence rate of the data-independent {$P$}-greedy algorithm in kernel-based approximation}, Dolomites Res. Notes Approx. \textbf{10} (2017), pp. 68--78.
	
	\bibitem{Santin24} 
	\textsc{G. Santin, T. Wenzel, B. Haasdonk}, 
	\emph{On the Optimality of Target-Data-Dependent Kernel Greedy Interpolation in Sobolev Reproducing Kernel Hilbert Spaces}, SIAM Journal on Numerical Analysis \textbf{62}(5) (2024), pp. 2249--2275.
	
	\bibitem{Schaback00} 
	\textsc{R. Schaback, H. Wendland}, 
	\emph{Adaptive greedy techniques for approximate solution of large {RBF} systems}, Numer. Algorithms \textbf{24}(3) (2000), pp. 239--254.
	
	\bibitem{Scheuerer11}
	\textsc{M. Scheuerer}, \emph{An alternative procedure for selecting a good value for the parameter c in RBF-interpolation}, Adv. Comput. Math. \textbf{34} (2011), pp. 105--126.
	
	\bibitem{Scheuerer13} 
	\textsc{M. Scheuerer, R. Schaback, M. Schlather}, 
	\emph{Interpolation of spatial data – A stochastic or a deterministic problem?}, European Journal of Applied Mathematics \textbf{24}(4) (2013), pp. 601--629.
	
	\bibitem{Srinivas10} 
	\textsc{N. Srinivas, A. Krause, S. Kakade, M. Seeger}, 
	\emph{Gaussian process optimization in the bandit setting: no regret and experimental design}, In International Conference on Machine Learning (2010), pp. 1015--1022.	
	
	\bibitem{Vakili22}
	\textsc{S. Vakili}, 
	\emph{Open Problem: Regret Bounds for Noise-Free Kernel-Based Bandits}, PMLR \textbf{178} (2022), pp. 5624--5629.
	
	\bibitem{venkatanath2015}
	\textsc{N. Venkatanath, D. Praneeth, B. Maruthi Chandrasekhar, S. S. Channappayya , S.M. Swarup},
	\emph{ Blind Image Quality Evaluation using perception based features.}, 2015 Twenty First National Conference on Communications (NCC), Mumbai, India (2015), pp. 1--6.
	
	
	
	\bibitem{Wendland05}
	\textsc{H. Wendland}, 
	\emph{Scattered Data Approximation}, 
	Cambridge Monogr. Appl. Comput. Math., vol. 17, Cambridge Univ. Press, Cambridge, 2005.
	
	\bibitem{Wendland10} 
	\textsc{H. Wendland}, 
	\emph{Multiscale analysis in Sobolev spaces on bounded domains}, Numer. Math. \textbf{116} (2010), pp. 493--517.
	
	\bibitem{Wendland05a} 
	\textsc{H. Wendland, C. Rieger}, 
	\emph{Approximate Interpolation with Applications
		to Selecting Smoothing Parameters}, Numer. Math. \textbf{101} (2005), pp. 729--748.
	
	\bibitem{Wenzel21} 
	\textsc{T. Wenzel, G. Santin, B. Haasdonk}, 
	\emph{A novel class of stabilized greedy kernel approximation algorithms: Convergence, stability and uniform point distribution}, J. Approx. Theory \textbf{262} (2021), 105508.	
	
	\bibitem{Wenzel23} 
	\textsc{T. Wenzel, G. Santin, B. Haasdonk}, 
	\emph{{A}nalysis of {T}arget {D}ata-{D}ependent {G}reedy {K}ernel {A}lgorithms: {C}onvergence {R}ates for $f$-, $f \cdot {P}$- and $f/{P}$-{G}reedy}, Constr. Approx. \textbf{57} (2023), pp. 45--74.	
	
	
	\bibitem{Wenzel24} 
	\textsc{T. Wenzel, F. Marchetti, E. Perracchione}, 
	\emph{Data-Driven Kernel Designs for Optimized Greedy Schemes: A Machine Learning Perspective}, SIAM Journal on Scientific Computing \textbf{46}(1) (2024), pp. C101--C126.	
	
	\bibitem{Wirtz13} 
	\textsc{D. Wirtz, B. Haasdonk}, 
	\emph{A vectorial kernel orthogonal greedy algorithm}, Dolomites Res. Approx. \textbf{6} (2013), pp. 83--100.	
	
	\bibitem{Wirtz15} 
	\textsc{D. Wirtz, N. Karajan, B. Haasdonk}, 
	\emph{Surrogate modeling of multiscale models using kernel methods}, Int. J. Numer. Methods Eng. \textbf{101}(1) (2015), pp. 1--28.	
	
	\bibitem{Zanzonico04} 
	\textsc{P.B. Zanzonico}, 
	\emph{Principles of nuclear medicine imaging: Planar, SPECT, PET, and hybrid modalities}, Radiographics \textbf{24}(6) (2004), pp. 1621--1641.	
	
	\bibitem{Zeng10}
	\textsc{G.L. Zeng}, 
	\emph{Medical Image Reconstruction: A Conceptual Tutorial}, Springer, 2010.
	
	
	\bibitem{Zorz24} 
	\textsc{A. Zorz et al.}, 
	\emph{Performance evaluation of the 3D-ring cadmium–zinc–telluride (CZT) StarGuide system according to the NEMA NU 1-2018 standard}, EJNMMI Physics \textbf{11} (2024), 69.	

\end{thebibliography}

\end{document}